\documentclass[11pt]{article}
\usepackage{amsmath,amsthm,amssymb}
\usepackage[T1]{fontenc}
\usepackage[utf8]{inputenc}
\usepackage{graphicx}
\usepackage{subcaption}
\graphicspath{{figures/}{empirical_path_law_AW/figures/}}
\usepackage{color}
\usepackage{epstopdf}
\usepackage{mathtools}
\mathtoolsset{showonlyrefs=true}
\usepackage{tikz}
\usepackage{enumerate}
\usepackage{enumitem}
\usepackage{algorithm}
\usepackage{algpseudocode}
\usepackage{mathrsfs}
\usepackage[normalem]{ulem}
\usepackage[authoryear,round]{natbib}
\usepackage[colorlinks=true,citecolor=db,linkcolor=db,urlcolor=blue,pdfstartview=FitH]{hyperref}
\usepackage{comment}
\usepackage{array}

\definecolor{db}{RGB}{0, 0, 130}
\definecolor{rp}{rgb}{0.25, 0, 0.75}
\definecolor{dg}{rgb}{0, 0.6, 0}

\newtheorem{theorem}{Theorem}[section]

\newtheorem{definition}{Definition}[section]

\newtheorem{corollary}[definition]{Corollary}

\newtheorem{assumption}[definition]{Assumption}
\newtheorem{lemma}[definition]{Lemma}

\newtheorem{proposition}[definition]{Proposition}
\newtheorem{remark}[definition]{Remark}

\def\Ac{\mathcal A}

\def\Cc{\mathcal C}
\def\Fc{\mathcal F}

\def\Lc{\mathcal L}

\def\Pc{\mathcal P}

\def\E{\mathbb{E}}

\def\N{\mathbb N}
\def\P{\mathbb P}
\def\R{\mathbb R}

\def\1{\mathbf{1}}
\newcommand{\cW}{\mathcal W}
\usepackage{booktabs}

\numberwithin{equation}{section}
\setlist{leftmargin=*,itemsep=2pt,topsep=4pt}
\def\Pp{\P}
\def\W{\mathbb{W}}
\def\AW{\mathcal{AW}}
\def\Can{\operatorname{Can}}
\def\Law{\operatorname{Law}}
\def\supp{\operatorname{supp}}

\def\TV{\operatorname{TV}}
\def\cF{\Fc}
\def\cA{\Ac}
\def\cC{\Cc}
\def\cP{\Pc}

\allowdisplaybreaks[1]
\hypersetup{pdftitle={Approximating It\^o Processes in Adapted Wasserstein Distance with Finite Support Laws},pdfauthor={Xihao He and Fengyi Yuan}}

\title{Approximating It\^o Processes in Adapted Wasserstein Distance with Finite Support Laws}
\author{
Xihao He\thanks{Department of Mathematics, University of Southern California. Email: \href{mailto:xihaohe@usc.edu}{\texttt{xihaohe@usc.edu}}.}
\and
Fengyi Yuan\thanks{Corresponding author. School of Science and Engineering, the Chinese University of Hong Kong (Shenzhen). Email: \href{mailto:yuanfengyi@cuhk.edu.cn}{\texttt{yuanfengyi@cuhk.edu.cn}}. Fengyi Yuan is supported by the Chinese University of Hong Kong (Shenzhen) start-up fund under UDF01004253.}%
}
\date{}

\begin{document}
\maketitle

\begin{abstract}
Motivated by finite-scenario surrogates in stopping and control problems, we study the approximation of It\^o diffusion laws by laws supported on $M$ paths in an adapted Wasserstein distance, which we term {\it adapted functional quantization}. Under the strict adapted Wasserstein distance, no finite-support approximation can achieve vanishing error, while even under its relaxed counterpart, naive empirical measures still fail to converge.
To resolve this, we develop a recursive conditional sampling algorithm and derive constructive upper bounds on the adapted quantization error.
In one dimension, this yields the Brownian rate $(\log M)^{-1/2}$ for every diffusion with continuous coefficients of linear growth and pathwise uniqueness.
In several dimensions, a quantitative Osgood stability condition with exponent $\beta$ yields the rate $(\log M)^{-\eta_p}$, where $\eta_p=\min\{1/2,\beta/(\beta+4p-2)\}$.
Without the Osgood condition, we construct a two-dimensional diffusion with smooth bounded coefficients whose adapted quantization error decays more slowly than $(\log M)^{-\gamma}$ for any $\gamma>0$.
\end{abstract}

\tableofcontents

\section{Introduction}

In dynamic stochastic control problems, ranging from optimal stopping and American option pricing to dynamic portfolio selection, a standard computational paradigm is to replace a complex continuous-time stochastic process with a tractable finite-scenario surrogate \citep{Pflug2001,BallyPages2003,BallyPagesPrintems2005,LohndorfWozabal2023,Pflug2023,BayraktarFengZhang2024,BayerPelizzariSchoenmakers2025,AghapourBayraktarYuan2025}.
When the underlying dynamics follow an It\^o diffusion $X=(X_t)_{t\in[0,1]}$ on $\cC_d:=C([0,1];\R^d)$, this amounts to approximating the path law $\Lc(X)$ by a discrete probability measure $\nu$ supported on at most $M$ trajectories.

Classically, the best $M$-point approximation of a probability distribution in an $L^p$-type metric is measured by $$e_{M,p}(X) :=\inf_{|\mathcal A|\leq M}\big(\E[d(X,\mathcal A)^p]\big)^{1/p}=\inf_{|\supp\nu|\leq M}\cW_p(\Lc(X),\nu),$$ whose asymptotic behavior is the subject of quantization \citep{Pollard1981,GrayNeuhoff1998,GrafLuschgy2000} or its infinite-dimensional counterpart, functional quantization \citep{LuschgyPages2002,DereichScheutzow2006}.
When $X$ is a process with continuous paths and $d$ is the uniform distance, the quantization error $e_{M,p}(X)$ is now well understood. For a broad class of It\^o processes, $e_{M,p}(X)$ decays at the universal Brownian rate $(\log M)^{-1/2}$, and this rate can be achieved constructively by the empirical measure of $M$ independent sample paths \citep{HeYuan2026}; see also \citet{DereichScheutzowSchottstedt2013,BoissardLeGouic2014,FournierGuillin2015,WeedBach2019,Lei2020} for empirical Wasserstein convergence in finite- and infinite-dimensional spaces.

However, it is a fundamental insight of modern stochastic control that the classical Wasserstein topology is too coarse for dynamic decision-making.
Two processes can be arbitrarily close in $\cW_p$ while possessing incompatible filtrations, leading to drastic discrepancies in optimal stopping values and hedging strategies \citep{PflugPichler2012,BackhoffEtAl2020Stability,BackhoffEtAl2020Topologies,AcciaioBackhoffZalashko2020}.
To ensure the continuity of value functions in stopping and control problems, the approximation must respect the flow of information.
This requirement is captured by bicausal optimal transport and the resulting adapted Wasserstein distance \citep{PflugPichler2012,BackhoffEtAl2020Stability,BackhoffEtAl2020Topologies,BBPSZ2025}.
Motivated by these developments, we study the approximation of It\^o diffusion laws in adapted Wasserstein distance by laws supported on at most $M$ paths, which we call {\it adapted functional quantization}. More specifically, we consider
\begin{equation}\label{intro:rate}
e^{\mathrm{ad}}_{M,p}(\mathbb X):=\inf_{|\supp\nu|\leq M}\AW_p(\mathbb X,\Can(\nu)),
\end{equation}
where $\mathbb X$ and $\Can(\nu)$ are canonical liftings of $X$ and $\nu$ needed to define the adapted distance. Importantly, since $e^{\mathrm{ad}}_{M,p}(\mathbb X)$ merely represents the theoretical best approximation under a given support budget, we aim to bound it {\it constructively} by designing an implementable sampling scheme with an explicit error decay rate, from which upper bounds on $e^{\mathrm{ad}}_{M,p}(\mathbb X)$ follow immediately.

\paragraph{Two foundational obstructions.}
Passing from classical to adapted functional quantization in continuous time immediately encounters two obstacles that have no counterpart under $\cW_p$.

First, the {\it strict} adapted Wasserstein distance $\AW_p^{(s)}$, which enforces bicausality without time delay at every instant $t\in[0,1]$, precludes any vanishing quantization error.
Whereas the natural filtration of a diffusion evolves continuously, the canonical filtration of a finitely supported law $\nu$ is piecewise constant and changes at most $|\supp\nu|-1$ times.
Building on \citet{BBPSZ2025}, Proposition~\ref{prop:strict-finite-support} shows that the only strictly bicausal coupling between $d$-dimensional Wiener measure $\W_d$ and any finite-support law $\nu$ is the independent product coupling.
Consequently, the strict adapted quantization error of Brownian motion remains bounded away from zero for all $M\in\N$.
To obtain a meaningful theory of finite-support approximation, one must allow a small temporal relaxation in the causality constraint.
Therefore, in \eqref{intro:rate} we take $\AW_p$ to be the {\it relaxed} adapted Wasserstein distance introduced by \citet{BBPSZ2025}, which allows a causality delay $\varepsilon\geq 0$ and penalizes it in the cost.

Second, even under the relaxed distance $\AW_p$, {raw empirical path sampling fails}.
If $\mu_M$ is the empirical measure of $M$ independent Brownian motions, almost surely the $M$ sample trajectories are pairwise distinct at time $1/n$ for every $n\in\N$.
Observing any infinitesimal initial segment of a trajectory under $\mu_M$ immediately reveals which of the $M$ paths has been selected, collapsing the canonical filtration at time $t=0+$ to the full terminal $\sigma$-field.
As a result, Proposition~\ref{prop:sampling-sharp} proves that $\AW_p(\W_d,\Can(\mu_M))$ stays bounded below by a positive constant for all $M\in\N$.
This continuous-time phenomenon parallels the well-known failure of raw empirical measures in discrete-time adapted estimation, where smoothing or clustering modifications are required to restore convergence \citep{BackhoffEtAl2022,AcciaioHou2024,LarssonParkWiesel2025,MirmominovWiesel2025}.

\paragraph{Main results.}
The aforementioned obstructions motivate us to develop and analyze a recursive conditional sampling scheme (Algorithm~\ref{alg:recursive-sampling}), through which we obtain constructive upper bounds for $e^{\mathrm{ad}}_{M,p}(\mathbb X)$.
This algorithm recursively samples new branches at each endpoint of the previously sampled paths.
Assuming exact conditional sampling,\footnote{Approximation errors from estimated or learned samplers, or from numerical discretization of path law sampling, require additional analysis and are left for future study.} Theorem~\ref{thm:recursive-sampling} derives the mean relaxed adapted Wasserstein error rates of such a sampled law, yielding the quantization upper bounds in Theorem~\ref{thm:quantization}.
In one dimension, these have the Brownian rate $(\log M)^{-1/2}$ under continuous coefficients of linear growth and strong well-posedness.
In several dimensions, the logarithmic upper rates are determined by the finite-time stability of a joint Osgood flow. We believe that this intuitive recursive sampling scheme offers a convenient alternative to the {\it adapted empirical measure} \citep[cf.][]{BackhoffEtAl2022}, whose direct extension to continuous time appears considerably more involved. Our proof constructs relaxed bicausal couplings interval by interval through local matching of path and endpoint labels, and controls accumulated endpoint perturbations using the scalar comparison principle or quantitative Osgood stability.

We complement these sampling guarantees with lower bounds that apply to every finite-support law.
Our lower-bound method uses the complexity of conditional laws of future observations.
Proposition~\ref{prop:counterexample} constructs Osgood diffusions with an explicit lower bound for every prescribed stability exponent $\beta>0$.
As $\beta\downarrow0$, the upper and lower logarithmic exponents agree to first order, establishing first-order sharpness of our adapted $L^1$ quantization bounds (Remark~\ref{rem:asymptotic-optimality}).
The same family contains examples with adapted $L^1$ quantization rates strictly slower than their classical counterparts (Remark~\ref{rem:classical-separation}), quantifying the additional approximation cost imposed by the flow of information.
Finally, Proposition~\ref{prop:slow-quantization} constructs a two-dimensional diffusion satisfying Assumption~\ref{ass:coefficients} alone whose adapted $L^1$ quantization error decays more slowly than any negative power of $\log M$.
This example demonstrates the essential role of quantitative stability in several dimensions. See also Table \ref{tab:summary} for a summary of our main results.

\begin{table}[htbp]
\centering
\caption{Summary of main results}
\label{tab:summary}
\small
\setlength{\tabcolsep}{4pt}
\renewcommand{\arraystretch}{1.12}
\begin{tabular*}{\linewidth}{@{\extracolsep{\fill}}lll@{}}
\toprule
Setting & Upper bound & Lower bound\\
\midrule
{\it Strict} $\AW$
 & -- & non-vanishing\enspace
    {\footnotesize (Prop.~\ref{prop:strict-finite-support})}\\
\addlinespace[8pt]
\multicolumn{3}{@{}l}{\emph{Relaxed} $\AW$}\\
\addlinespace[3pt]
\quad Raw empirical
 & -- & non-vanishing\enspace
    {\footnotesize (Prop.~\ref{prop:sampling-sharp})}\\
\addlinespace[6pt]
\quad $d=1$
 & $(\log M)^{-1/2}$\enspace
     {\footnotesize (Thm.~\ref{thm:quantization}\textnormal{(i)})}
 & --\\
\addlinespace[6pt]
\quad $d\geq2$, $\beta$-Osgood
 & $(\log M)^{-\eta_p}$\,\textsuperscript{\dag}\enspace
     {\footnotesize (Thm.~\ref{thm:quantization}\textnormal{(ii)})}
 & $(\log M)^{-\beta/2-\beta^2/2}$\enspace
     {\footnotesize ($p=1$, Prop.~\ref{prop:counterexample})}\\
\addlinespace[6pt]
\quad $d\geq2$, no Osgood
 & --
 & $\gg(\log M)^{-\gamma}$, $\forall\gamma>0$\enspace
     {\footnotesize (Prop.~\ref{prop:slow-quantization})}\\
\bottomrule
\end{tabular*}
\par\vspace{3pt}
\parbox{\linewidth}{\footnotesize
  \textsuperscript{\dag}$\eta_p=\min\{1/2,\beta/(\beta+4p-2)\}$.
  Both upper bounds are attained by recursive sampling.}
\end{table}

\paragraph{Contributions.}
In summary, this paper initiates the study of functional quantization under adapted Wasserstein distances and makes the following contributions:
\begin{enumerate}[label=\textnormal{(\roman*)}]
\item 
We propose recursive conditional sampling and establish the first quantitative adapted Wasserstein convergence guarantee for a finite-support approximation of a continuous-time diffusion, with explicit error bounds in the depth and branching parameters.
\item 
We recover the Brownian rate $(\log M)^{-1/2}$ for all well-posed scalar diffusions, and obtain the explicit rate $(\log M)^{-\eta_p}$ in higher dimensions under a quantitative Osgood condition, with explicit lower bounds that are matching to first order for $p=1$ as the stability exponent vanishes.
These results reveal that, unlike its classical counterpart, the adapted quantization rate of a diffusion is governed by the quantitative stability of its flow.
\item 
We exhibit examples in which adapted quantization is provably slower than classical quantization, and show that in the absence of flow stability the adapted rate can be slower than any negative power of $\log M$ already in dimension two.
This quantifies, for the first time in continuous time, the additional approximation cost imposed by the flow of information.
\end{enumerate}

\paragraph{Organization of the paper.}
Section~\ref{sec:notation} introduces the strict and relaxed adapted Wasserstein distances, formulates our standing assumptions, and proves the obstruction to strict finite-support quantization (Proposition~\ref{prop:strict-finite-support}).
Section~\ref{sec:main} states our main results: the adapted quantization upper bounds (Theorem~\ref{thm:quantization}), the failure of naive empirical measures (Proposition~\ref{prop:sampling-sharp}), the recursive conditional sampling algorithm and its mean error bounds (Theorem~\ref{thm:recursive-sampling}), and the lower bound constructions establishing sharpness and slow rates (Propositions~\ref{prop:counterexample} and~\ref{prop:slow-quantization}).
Section~\ref{sec:proofs} outlines the overarching proof strategies and provides the complete proofs of all main and auxiliary results.

\section{Problem Formulation}
\label{sec:notation}

\subsection{Strict and relaxed adapted Wasserstein distances}
For $d\in\N$, let $\cC_d:=C([0,1];\R^d)$, equipped with its Borel $\sigma$-field and the uniform norm $\|x\|_\infty:=\sup_{0\leq t\leq1}|x_t|$.
We use $|\cdot|$ for the Euclidean norm and $\|\cdot\|_{\mathrm{HS}}$ for the Hilbert--Schmidt norm.
Let $\cP_p(\cC_d)$ denote the path laws with a finite $p$-th uniform-norm moment and $\Lc(X)$ the law of $X$.
Write $\cW_p$ for the ordinary Wasserstein distance on path space, $\gamma_d$ for $d$-dimensional Wiener measure, and $\W_d:=\Can(\gamma_d)$.

For $\mu\in\cP_p(\cC_d)$, let $\Can(\mu)$ be the coordinate process with law $\mu$ and its completed, right-continuous natural filtration.
With $\mathcal N_\mu$ denoting the subsets of $\mu$-null sets,
\begin{equation}
 \cF_t^\mu :=\bigcap_{t<u\leq1}\bigl(\sigma(X_s:0\leq s\leq u)\vee\mathcal N_\mu\bigr),\quad 0\leq t<1,\qquad
 \cF_1^\mu:=\sigma(X_s:0\leq s\leq1)\vee\mathcal N_\mu.
\end{equation}
On a coupling space, marginal filtrations are identified with their pullbacks and completed modulo the coupling's null sets.

\begin{definition}\label{def:bicausal-couplings}
A coupling $\pi$ of filtered processes $\mathbb X$ and $\mathbb Y$ is $\varepsilon$-causal from $\mathbb X$ to $\mathbb Y$ if, for every $t\in[0,1]$,
the $\sigma$-fields $\cF_t^Y$ and $\cF_1^X$ are conditionally independent given $\cF_{(t+\varepsilon)\wedge1}^X$.
It is $\varepsilon$-bicausal if this condition also holds with the marginals interchanged.
Write $\Pi_{bc}^{\varepsilon}(\mathbb X,\mathbb Y)$ for the set of such bicausal couplings.
\end{definition}

For $p\in[1,\infty)$, define the strict and relaxed adapted Wasserstein distances by
\begin{align}
 & \AW_p^{(s)}(\mathbb X,\mathbb Y)
 :=\inf_{\pi\in\Pi_{bc}^0}(\E_\pi\|X-Y\|_\infty^p)^{1/p},\notag\\
 & \AW_p(\mathbb X,\mathbb Y)
 :=\inf_{\varepsilon\geq0}\left\{\varepsilon+
 \inf_{\pi\in\Pi_{bc}^{\varepsilon}}(\E_\pi\|X-Y\|_\infty^p)^{1/p}\right\}.
 \label{eq:aw-definition}
\end{align}
Here $\Pi_{bc}^{\varepsilon}:=\Pi_{bc}^{\varepsilon}(\mathbb X,\mathbb Y)$.
In particular,
\begin{equation}
 \cW_p(\mu,\nu)\leq\AW_p(\Can(\mu),\Can(\nu))\leq\AW_p^{(s)}(\Can(\mu),\Can(\nu)).
\end{equation}
The adapted and ordinary $M$-point quantization errors are
\begin{equation}\label{eq:quantization-errors}
 e^{\mathrm{ad}}_{M,p}(\mathbb X):=\inf_{|\supp\nu|\leq M}\AW_p(\mathbb X,\Can(\nu)),\quad
 e_{M,p}(X):=\inf_{1\leq|A|\leq M}\left(\E\min_{a\in A}\|X-a\|_\infty^p\right)^{1/p}.
\end{equation}
The infima run over $\nu\in\cP_p(\cC_d)$ and $A\subset\cC_d$, respectively.
Thus $e^{\mathrm{ad}}_{M,p}(\Can(\Lc(X)))\geq e_{M,p}(X)$.

\subsection{Assumptions}
\label{sec:assumptions}
Let $B$ be an $m$-dimensional Brownian motion and consider a diffusion process $X = (X_s)_{s\in[0,1]}$ that follows
\begin{equation}\label{eq:sde}
 X_t=x_0+\int_0^t b(s,X_s)\,\mathrm ds+\int_0^t\sigma(s,X_s)\,\mathrm dB_s,\qquad 0\leq t\leq1,
\end{equation}
with fixed dimensions $d,m\in\N$ satisfying the following assumptions:

\begin{assumption}\label{ass:coefficients}
The initial state $x_0\in\R^d$ is deterministic.
The coefficients $b:[0,1]\times\R^d\to\R^d$ and $\sigma:[0,1]\times\R^d\to\R^{d\times m}$ are continuous and satisfy, for some $K<\infty$,
\begin{equation}
 |b(t,x)|+\|\sigma(t,x)\|_{\mathrm{HS}}\leq K(1+|x|),\qquad (t,x)\in[0,1]\times\R^d.
\end{equation}
For every deterministic starting pair $(s,x)\in[0,1]\times\R^d$, the equation on $[s,1]$ with $X_s=x$ admits a nonexplosive strong solution $X^{s,x}$ and pathwise uniqueness holds.
\end{assumption}
Set $\mathbb X:=\Can(\Lc(X))$.
Uniform linear growth gives $\E\|X\|_\infty^q<\infty$ for every finite $q\geq1$.

\begin{assumption}\label{ass:osgood}
There are nonnegative continuous functions $\ell,k:[0,1]\to[0,\infty)$ and an increasing concave function $\rho:[0,\infty)\to[0,\infty)$ with $\rho(0)=0$, $\rho(u)>0$ for $u>0$, and $\int_{0+}\mathrm du/\rho(u)=\infty$, such that
\begin{equation}\label{eq:joint-osgood}
 2\langle x-y,b(t,x)-b(t,y)\rangle+\|\sigma(t,x)-\sigma(t,y)\|_{\mathrm{HS}}^2
 \leq\ell(t)|x-y|^2+k(t)\rho(|x-y|^2)
\end{equation}
for all $t\in[0,1]$ and $x,y\in\R^d$.
\end{assumption}
Define the Osgood flow by
\begin{align}
 & z'(t;a)=\ell(t)z(t;a)+k(t)\rho(z(t;a)),\qquad z(0;a)=a>0,
 \label{eq:comparison-flow}\\
 & \Phi(a):=z(1;a)\quad\text{for }a>0, \qquad \Phi(0):=0.
\end{align}
Concavity bounds the vector field by a function of linear growth, and the Osgood condition gives $\Phi(a)\to0$ as $a\downarrow0$.
We quantify its order near zero as follows.

\begin{assumption}\label{ass:flow-rate}
For some $\beta>0$, the Osgood flow satisfies
\begin{equation}\label{eq:flow-condition}
 \Phi(a)=O\bigl([\log(1/a)]^{-\beta}\bigr),\qquad a\downarrow0.
\end{equation}
\end{assumption}

We next give examples of moduli in Assumption~\ref{ass:osgood} for which the Osgood flow satisfies \eqref{eq:flow-condition} for every $\beta>0$.
The following proposition covers moduli dominated by any fixed finite product of iterated logarithms.

\begin{proposition}\label{prop:iterated-log-flow}
Let $\ell,k:[0,1]\to[0,\infty)$ be continuous, and let $\rho:[0,\infty)\to[0,\infty)$ be increasing and concave, with $\rho(0)=0$ and $\rho(u)>0$ for $u>0$.
Suppose that, for some constant $C>0$, some fixed integer $j\geq1$, and all sufficiently small $u>0$,
\begin{equation}\label{eq:iterated-log-modulus}
 \rho(u)\leq Cu\prod_{i=1}^j\log^{(i)}(1/u),
\end{equation}
where $\log^{(i)}$ denotes the $i$-fold iterated logarithm.
Then $\int_{0+}\mathrm du/\rho(u)=\infty$, and the Osgood flow \eqref{eq:comparison-flow} satisfies \eqref{eq:flow-condition} for every $\beta>0$.
In particular, these conclusions hold for $\rho(u)=u$.
\end{proposition}

\begin{proof}
See Section~\ref{sec:notation-proofs}.
\end{proof}

\subsection{A negative result under strict adapted Wasserstein distance}
\label{sec:strict-quantization}

We first explain why temporal relaxation is needed for finite-support approximation.
For Wiener measure, the following consequence of \citet[Proposition~C.2]{BBPSZ2025} shows that the strict quantization error remains positive, regardless of the number of support paths.

\begin{proposition}\label{prop:strict-finite-support}
Let $d\in\N$, $p\in[1,\infty)$, and let $\nu\in\cP_p(\cC_d)$ have finite support.
Every strictly bicausal coupling between $\W_d$ and $\Can(\nu)$ is the product coupling $\gamma_d\otimes\nu$.
Consequently, if $B$ is a standard $d$-dimensional Brownian motion, then
\begin{equation}\label{eq:strict-finite-support-cost}
 \AW_p^{(s)}\bigl(\W_d,\Can(\nu)\bigr)^p
 =\int_{\cC_d}\E\|B-y\|_\infty^p\,\nu(\mathrm dy)
 \geq\E\|B\|_\infty^p,
\end{equation}
and, for every $M\in\N$,
\begin{equation}\label{eq:strict-quantization-obstruction}
 \inf_{\substack{\nu\in\cP_p(\cC_d)\\|\supp\nu|\leq M}}
 \AW_p^{(s)}\bigl(\W_d,\Can(\nu)\bigr)
 =\bigl(\E\|B\|_\infty^p\bigr)^{1/p}>0.
\end{equation}
\end{proposition}

\begin{proof}
See Section~\ref{sec:notation-proofs}.
\end{proof}

Our main results give explicit vanishing quantization bounds in the relaxed adapted Wasserstein distance.

\section{Main results}
\label{sec:main}

\subsection{Quantization upper bounds}
\label{sec:quantization-rates}
We first state the adapted functional quantization upper bounds as a part of our main results.
Section~\ref{sec:recursive-sampling} presents further details behind the proof of these results. In particular, it gives a recursive conditional sampling construction that attains these rates in mean relaxed adapted Wasserstein error.

\begin{theorem}
\label{thm:quantization}
Suppose Assumption~\ref{ass:coefficients} holds and fix $p\in[1,\infty)$.
There is a constant $C_p>0$ such that the following bounds hold for every sufficiently large integer $M$.
\begin{enumerate}[label=\textnormal{(\roman*)}]
\item\label{item:scalar-quantization} If $d=1$, then
\begin{equation}\label{eq:scalar-quantization-consequence}
 e^{\mathrm{ad}}_{M,p}(\mathbb X)\leq C_p(\log M)^{-1/2}.
\end{equation}
\item\label{item:multidimensional-quantization} If $d\geq2$ and Assumptions~\ref{ass:osgood} and~\ref{ass:flow-rate} hold, then
\begin{equation}\label{eq:main-multidimensional-quantization}
 e^{\mathrm{ad}}_{M,p}(\mathbb X)\leq C_p(\log M)^{-\eta_p},\qquad
 \eta_p:=\min\left\{\frac12,\frac{\beta}{\beta+4p-2}\right\}.
\end{equation}
\end{enumerate}
The constant and threshold may depend on the coefficients, the dimensions, $x_0$, $p$, and the constants in the applicable assumptions, but are independent of $M$.
\end{theorem}

\begin{proof}
These bounds follow from Theorem~\ref{thm:recursive-sampling} by choosing the sampling parameters as described at the end of Section~\ref{sec:recursive-sampling}.
\end{proof}

\begin{corollary}\label{cor:iterated-log-rate}
Suppose Assumptions~\ref{ass:coefficients} and~\ref{ass:osgood} hold with a modulus covered by Proposition~\ref{prop:iterated-log-flow}.
Then, for every $p\in[1,\infty)$, there is a constant $C_p<\infty$ such that, for all sufficiently large $M$,
\begin{equation}\label{eq:iterated-log-rate}
 e^{\mathrm{ad}}_{M,p}(\mathbb X)\leq C_p(\log M)^{-1/2}.
\end{equation}
\end{corollary}

\begin{proof}
Proposition~\ref{prop:iterated-log-flow} gives Assumption~\ref{ass:flow-rate} for every $\beta>0$.
Choose $\beta\geq4p-2$ and apply Theorem~\ref{thm:quantization}\textnormal{(ii)} to obtain \eqref{eq:iterated-log-rate} when $d\geq2$.
Part~\textnormal{(i)} gives the same bound for $d=1$.
\end{proof}

\begin{corollary}[Wiener measure]\label{cor:wiener}
For every $d\in\N$ and $p\in[1,\infty)$, there are constants $0<c_{p,d}\leq C_{p,d}<\infty$ such that, for all sufficiently large $M$,
\begin{equation}
 c_{p,d}(\log M)^{-1/2}
 \leq e^{\mathrm{ad}}_{M,p}(\W_d)
 \leq C_{p,d}(\log M)^{-1/2}.
\end{equation}
\end{corollary}

\begin{proof}
The upper bound follows from Corollary~\ref{cor:iterated-log-rate} with $\rho(u)=u$.
Let $B$ be a standard $d$-dimensional Brownian motion and $B^1$ its first coordinate.
Coordinate projection and the classical supremum-norm quantization bound in \citet[Theorem~1.1]{DereichScheutzow2006}, with Hurst parameter $H=1/2$, give
$
 e^{\mathrm{ad}}_{M,p}(\W_d)\geq e_{M,p}(B)
 \geq e_{M,p}(B^1)\geq c_p(\log M)^{-1/2}.
$
\end{proof}

\begin{remark}
Independent scalar equations satisfying Assumption~\ref{ass:coefficients} also have joint adapted quantization rate $C_{p,d}(\log M)^{-1/2}$.
Apply the scalar quantization conclusion \eqref{eq:scalar-quantization-consequence} coordinatewise, take product couplings on a common grid, and assign at most $\lfloor M^{1/d}\rfloor$ paths to each coordinate.
\end{remark}

\begin{remark}\label{rem:brownian-scope}
Our results show that the Brownian quantization rate $(\log M)^{-1/2}$ extends well beyond Lipschitz coefficients.
In one dimension, Theorem~\ref{thm:quantization}\textnormal{(i)} imposes no prescribed modulus of continuity beyond Assumption~\ref{ass:coefficients}.
It covers, for example, the square-root short-rate model of \citet{CoxIngersollRoss1985}, with a nonnegative initial state and the square-root coefficient extended to $\R$ by $\sqrt{x^+}$, where $x^+:=\max\{x,0\}$.

In several dimensions, Assumption~\ref{ass:osgood} controls the drift through a one-sided inequality and permits moduli larger than linear near zero.
Related non-Lipschitz conditions for strong well-posedness are studied by \citet{FangZhang2005}.
For overdamped Langevin dynamics $\mathrm dX_t=-\nabla V(X_t)\,\mathrm dt+\sqrt2\,\mathrm dB_t$, convexity of $V$ makes the drift contribution in \eqref{eq:joint-osgood} nonpositive.
Corollary~\ref{cor:iterated-log-rate} therefore applies with $\rho(u)=u$ whenever $V\in C^1$ is convex and $\nabla V$ has at most linear growth, even if $\nabla V$ is not Lipschitz.
An example is $V(x)=q^{-1}\sum_{i=1}^d|x_i|^q$ with $q\in(1,2)$.
\end{remark}

\subsection{Recursive sampling and mean approximation error}
\label{sec:recursive-sampling}
Quantization describes the smallest approximation error achievable by a law supported on at most $M$ paths.
Therefore, it is important to seek an explicit and implementable construction that attains the quantization upper bounds.
For the ordinary Wasserstein distance, independent path sampling usually suffices for this purpose (in mean Wasserstein error); see \citep[Theorems~2.1 and~2.2]{HeYuan2026}.
However, this empirical construction fails in the relaxed adapted distance, even for Brownian motion.

\begin{proposition}\label{prop:sampling-sharp}
Let $d\in\N$, let $(B^i)_{i\geq1}$ be independent standard $d$-dimensional Brownian motions, and set
\begin{equation}
 \mu_M:=\frac1M\sum_{i=1}^M\delta_{B^i},\qquad M\in\N.
\end{equation}
On an event of probability one, for every $M\in\N$ and every $p\in[1,\infty)$,
\begin{equation}\label{eq:raw-empirical-lower}
 \AW_p\bigl(\W_d,\Can(\mu_M)\bigr)
 \geq\min\{1,m_{p,d}\}>0,
 \qquad m_{p,d}:=\bigl(\E|G|^p\bigr)^{1/p},
\end{equation}
where $G\sim\mathcal N(0,I_d)$.
In particular, the mean relaxed adapted Wasserstein error does not tend to zero.
\end{proposition}

\begin{proof}
See Section~\ref{sec:raw-empirical-proof}.
\end{proof}

This obstruction motivates us to study {\it recursive conditional sampling}: at each grid time, we extend every sampled prefix by several conditionally independent branches, producing a tree of paths.
Theorem~\ref{thm:recursive-sampling} bounds the mean relaxed adapted Wasserstein error of the resulting random law.
Optimizing the tree depth and branching number yields the quantization upper bounds in Theorem~\ref{thm:quantization}, as explained at the end of this subsection.

\paragraph{Recursive conditional sampling.}
Fix integers $n\geq 1$ and $L\geq 2$. We divide $[0,1]$ into $n$ intervals of
length $h:=1/n$, with grid points $t_k:=kh$ for $k=0,\ldots,n$. For $0\leq k<n$
and a starting point $y$, let
\begin{equation}
  P_k(y,\cdot):=\Lc\bigl((X_t^{t_k,y})_{t_k\leq t\leq t_{k+1}}\bigr)
\end{equation}
be the law of the exact diffusion path on $[t_k,t_{k+1}]$ that starts at $y$
at time $t_k$.

\emph{Idea of the construction.}
The construction has $n$ steps, one for each interval.
In the first step, we simulate $L$ independent paths on $[0,t_1]$, all starting
at $x_0$.
In step $k+1$, we take each path built so far; it is defined on $[0,t_k]$.
We look at its value at time $t_k$. From this value, we simulate $L$ new
independent paths on $[t_k,t_{k+1}]$, and we attach each of them to the old path.
We call a path that is simulated on a single interval $[t_k,t_{k+1}]$
a \emph{branch}.
Thus every path splits into $L$ paths at each grid point, and after $n$ steps we
obtain $L^n$ continuous paths on $[0,1]$.

\emph{Labels.}
We label the paths by finite sequences. Let $\mathcal I_0:=\{\varnothing\}$ and
$\mathcal I_k:=\{1,\ldots,L\}^k$ for $k\geq 1$. For
$\mathbf i=(i_1,\ldots,i_k)\in\mathcal I_k$ and $j\in\{1,\ldots,L\}$, we write
$\mathbf i j:=(i_1,\ldots,i_k,j)\in\mathcal I_{k+1}$.
For each $\mathbf i\in\mathcal I_k$, the path $g_{\mathbf i}$ is defined on
$[0,t_k]$. The label $\mathbf i$ records which branch was chosen in each of the
first $k$ steps.

\emph{Formal construction.}
We start with $g_{\varnothing}(0):=x_0$. For $k=0,\ldots,n-1$ and every
$\mathbf i\in\mathcal I_k$, we do the following.
\begin{enumerate}
  \item Draw $L$ branches $S_{\mathbf i,1},\ldots,S_{\mathbf i,L}$ from
        $P_k\bigl(g_{\mathbf i}(t_k),\cdot\bigr)$.
  \item For $j=1,\ldots,L$, define the path $g_{\mathbf i j}$ on $[0,t_{k+1}]$ by
        \begin{equation}
          g_{\mathbf i j}(t):=
          \begin{cases}
            g_{\mathbf i}(t), & t\in[0,t_k],\\
            S_{\mathbf i,j}(t), & t\in[t_k,t_{k+1}].
          \end{cases}
        \end{equation}
        We write this as $g_{\mathbf i j}=g_{\mathbf i}\star S_{\mathbf i,j}$
        (concatenation at time $t_k$).
\end{enumerate}
Each branch starts at the end point of the path it extends, that is,
$S_{\mathbf i,j}(t_k)=g_{\mathbf i}(t_k)$. So the two cases agree at $t=t_k$,
and $g_{\mathbf i j}$ is a well-defined continuous path.
The randomness is as follows: conditionally on all paths constructed in the
first $k$ steps, the new branches $S_{\mathbf i,j}$
($\mathbf i\in\mathcal I_k$, $j=1,\ldots,L$) are independent, and
$S_{\mathbf i,j}$ has law $P_k\bigl(g_{\mathbf i}(t_k),\cdot\bigr)$.

\begin{figure}[htbp]
\centering
\begin{subfigure}[t]{0.49\textwidth}
    \centering
    \includegraphics[width=\linewidth]{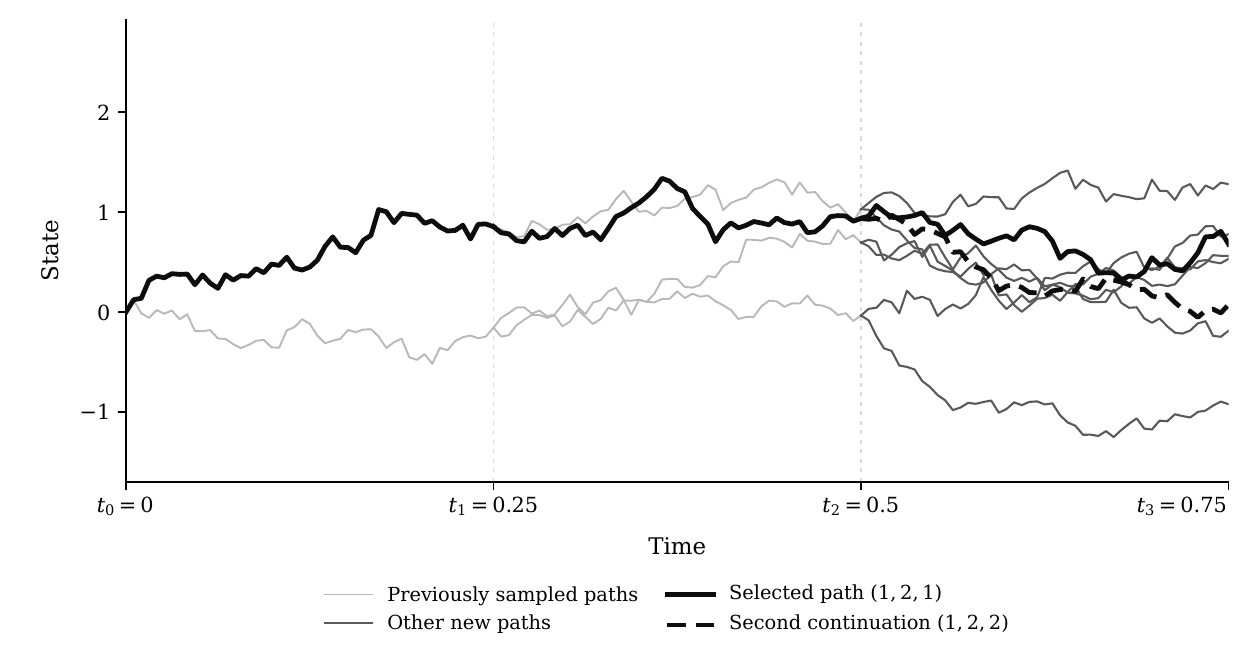}
    \caption{Each of four prefixes receives two branches. The bold prefix $(1,2)$ branches into $(1,2,1)$ (solid) and $(1,2,2)$ (dashed).}
\end{subfigure}\hfill
\begin{subfigure}[t]{0.49\textwidth}
    \centering
    \includegraphics[width=\linewidth]{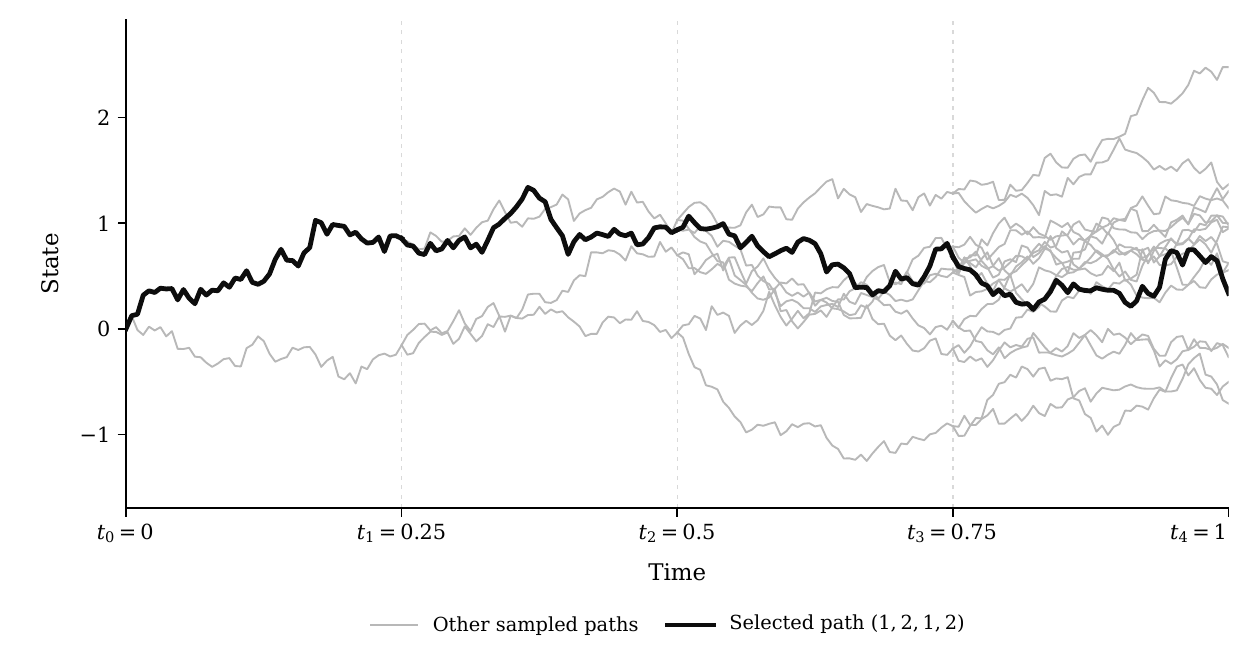}
    \caption{The $16$ indexed complete paths, each with mass $1/16$. The bold path is $g_{(1,2,1,2)}$.}
\end{subfigure}
\caption{Recursive conditional sampling with $n=4$ and $L=2$.}
\label{fig:recursive-sampling}
\end{figure}

After $n$ steps, we have $L^n$ continuous paths $g_{\mathbf i}$,
$\mathbf i\in\mathcal I_n$, on $[0,1]$. The random finite-support law is the
uniform measure on these paths,
\begin{equation}
  \mu^{n,L}:=\frac{1}{L^n}\sum_{\mathbf i\in\mathcal I_n}\delta_{g_{\mathbf i}}.
\end{equation}
Algorithm~\ref{alg:recursive-sampling} gives the pseudocode, and
Figure~\ref{fig:recursive-sampling} illustrates the construction.

\begin{algorithm}[htbp]
\caption{Recursive conditional sampling}
\label{alg:recursive-sampling}
\begin{algorithmic}[1]
\Require Integers $n\geq1$, $L\geq2$, initial state $x_0$, and conditional path samplers $P_0,\ldots,P_{n-1}$.
\State $g_{\varnothing}(0)\gets x_0$, $\mathcal P_0\gets[(\varnothing,g_{\varnothing})]$
\For{$k=0,\ldots,n-1$}
    \State $\mathcal P_{k+1}\gets[\,]$
    \ForAll{$(\mathbf i,g_{\mathbf i})\in\mathcal P_k$}
        \State Sample $S_{\mathbf i,1},\ldots,S_{\mathbf i,L}$ independently from $P_k(g_{\mathbf i}(t_k),\cdot)$
        \For{$j=1,\ldots,L$}
            \State Append $(\mathbf i j,g_{\mathbf i}\star S_{\mathbf i,j})$ to $\mathcal P_{k+1}$
        \EndFor
    \EndFor
\EndFor
\State \Return the uniform measure on the indexed paths in $\mathcal P_n$
\end{algorithmic}
\end{algorithm}

We measure the average approximation error by
\begin{equation}
 \mathcal E_{n,L,p}
 :=\E_{\mathcal S}[\AW_p\bigl(\mathbb X,\Can(\mu^{n,L})\bigr)],
\end{equation}
where $\E_{\mathcal S}$ denotes expectation with respect to the randomness of the recursive sampling procedure.

\begin{theorem}
\label{thm:recursive-sampling}
Suppose Assumption~\ref{ass:coefficients} holds and fix $p\in[1,\infty)$.
There exist $n_p\in\N$ and $a_p,C_p>0$ such that the following bounds hold for all integers $n\geq n_p$ and $L\geq\max\{2,n^{a_p}\}$.
\begin{enumerate}[label=\textnormal{(\roman*)}]
\item If $d=1$, then
\begin{equation}\label{eq:scalar-sampling-rate}
 \mathcal E_{n,L,p}\leq C_p\left[\frac1n+\frac1{\sqrt{n\log L}}\right].
\end{equation}
\item If $d\geq2$ and Assumptions~\ref{ass:osgood} and~\ref{ass:flow-rate} hold, then
\begin{equation}\label{eq:multidimensional-sampling-rate}
 \mathcal E_{n,L,p}\leq C_p\left[\frac1n+\frac1{\sqrt{n\log L}}
 +(\log L)^{-\beta/(4p-2)}\right].
\end{equation}
\end{enumerate}
The constants and thresholds may depend on the coefficients, the dimensions, $x_0$, $p$, and the constants in the applicable assumptions, but are independent of $n,L$.
\end{theorem}

\begin{proof}
See Section~\ref{sec:main-proofs}.
\end{proof}

\paragraph{Relation to quantization.}
Every realization of $\mu^{n,L}$ has support size at most $L^n$.
Consequently, whenever $L^n\leq M$, the definition of the quantization error gives
\begin{equation}
 e^{\mathrm{ad}}_{M,p}(\mathbb X)\leq\mathcal E_{n,L,p}.
\end{equation}
For every sufficiently large integer $M$, one can choose $n_M,L_M$ with $L_M^{n_M}\leq M$ such that, in case \textnormal{(i)} of Theorem~\ref{thm:recursive-sampling},
\begin{equation}
 e^{\mathrm{ad}}_{M,p}(\mathbb X)\leq\mathcal E_{n_M,L_M,p}
 \leq C_p(\log M)^{-1/2},
\end{equation}
and, in case \textnormal{(ii)},
\begin{equation}
 e^{\mathrm{ad}}_{M,p}(\mathbb X)\leq\mathcal E_{n_M,L_M,p}
 \leq C_p(\log M)^{-\eta_p},
\end{equation}
where $\eta_p$ is defined in Theorem~\ref{thm:quantization}.
Specifically, take $n_M$ of order $(\log M)^{1/2}$ in one dimension and $(\log M)^{\eta_p}$ in higher dimensions, with $L_M:=\lfloor M^{1/n_M}\rfloor$.
For sufficiently large $M$, these choices satisfy the conditions on $n,L$ in Theorem~\ref{thm:recursive-sampling} and the support budget $L_M^{n_M}\leq M$.
These estimates prove Theorem~\ref{thm:quantization} and show that the sampled laws attain its rates in mean relaxed adapted Wasserstein error.

\subsection{Lower bounds and separation from classical quantization}
\label{sec:osgood-optimality}
To derive the upper bounds, we combine classical path approximation with stability estimates for endpoint errors, using order preservation in one dimension and the Osgood flow in several dimensions.
Our lower bounds show that this stability is essential through a common mechanism.
A law supported on $M$ paths has at most $M$ conditional laws of the future at any fixed time.
If the filtration of $\mathbb X$ is constant on a short time interval, a bicausal coupling with a smaller delay cannot use future information there.
The adapted error is then at least the error of approximating the random conditional law of a Lipschitz observation of $\mathbb X$ by $M$ fixed laws (Lemma~\ref{lem:conditional-law-lower}).
Consequently, Corollary~\ref{cor:conditional-law-packing} gives a large adapted quantization error for a diffusion that realizes many well-separated conditional laws, each with non-negligible probability.
Stability limits the number of such laws by controlling how strongly the flow separates future laws arising from nearby present states.

The next result gives an explicit lower bound for each prescribed Osgood stability exponent $\beta>0$.

\begin{proposition}\label{prop:counterexample}
For every $\beta>0$, there exists a diffusion $X$ in dimension $d=\lceil4/\beta^2\rceil+2$, with bounded coefficients and $x_0=0$, satisfying Assumptions~\ref{ass:coefficients}, \ref{ass:osgood}, and \ref{ass:flow-rate} with exponent $\beta$, such that, for some constants $c_\beta,C_\beta>0$ and all sufficiently large $M$,
\begin{equation}\label{eq:counterexample-rate}
 c_\beta(\log M)^{-\beta/2-\beta^2/2}
 \leq e^{\mathrm{ad}}_{M,1}(\Can(\Lc(X)))
 \leq C_\beta(\log M)^{-\eta_1},
\end{equation}
where $\eta_1=\min\{1/2,\beta/(\beta+2)\}$.
The constants and the threshold for $M$ may depend on $\beta$.
\end{proposition}

\begin{proof}
See Section~\ref{sec:counterexample-proof}.
\end{proof}

\begin{remark}\label{rem:asymptotic-optimality}
Denote by $\alpha_\beta:=\beta/2+\beta^2/2$ the exponent in \eqref{eq:counterexample-rate}.
As $\beta\downarrow 0$, we have
\begin{equation}
 \frac{\alpha_\beta}{\eta_1}=1+\frac32\beta+\frac12\beta^2\longrightarrow1.
\end{equation}
Thus, at least for $p=1$, the exponent in Theorem~\ref{thm:quantization}\textnormal{(ii)} is asymptotically sharp to first order as $\beta\downarrow0$.
In particular, for every $\varepsilon>0$, choosing $\beta>0$ sufficiently small gives an example with lower bound $c_\beta(\log M)^{-\eta_1-\varepsilon}$.
\end{remark}

\begin{remark}\label{rem:classical-separation}
The diffusions in Proposition~\ref{prop:counterexample} have bounded coefficients, so \citet[Theorem~2.1]{HeYuan2026} yields the classical bound $e_{M,1}(X)\leq C(\log M)^{-1/2}$.
For $0<\beta<(\sqrt5-1)/2$, we have $\alpha_\beta<1/2$ and obtain
\begin{equation}
 \frac{e^{\mathrm{ad}}_{M,1}(\Can(\Lc(X)))}{e_{M,1}(X)}
 \geq c_\beta(\log M)^{1/2-\alpha_\beta}
 \longrightarrow\infty.
\end{equation}
Thus these diffusions have a strictly slower adapted quantization rate than their classical rate, showing that the information constraint imposes an additional approximation cost for weakly stable diffusions. See also the discussion in Remark~\ref{rem:stability-governs}.
\end{remark}

\subsection{A slow rate without Osgood stability}
\label{sec:slow-quantization}

In one dimension, Theorem~\ref{thm:quantization}\textnormal{(i)} yields the Brownian rate under Assumption~\ref{ass:coefficients} alone.
In the next result, we show that this conclusion fails in dimension two: without a quantitative stability condition, the adapted quantization error can decay more slowly than every negative power of $\log M$, even for smooth bounded coefficients.

\begin{proposition}\label{prop:slow-quantization}
There are bounded smooth coefficients $b:[0,1]\times\R^2\to\R^2$ and $\sigma:[0,1]\times\R^2\to\R^{2\times1}$ satisfying Assumption~\ref{ass:coefficients} with $x_0=0$, such that the corresponding solution $X$ satisfies
\begin{equation}\label{eq:slow-quantization-rate}
 \lim_{M\to\infty}(\log M)^\beta
 e^{\mathrm{ad}}_{M,1}(\Can(\Lc(X)))=+\infty,
 \qquad \forall \beta>0.
\end{equation}
\end{proposition}

\begin{proof}
See Section~\ref{sec:slow-quantization-proof}.
\end{proof}

\begin{remark}\label{rem:stability-governs}
Theorem~\ref{thm:quantization}\textnormal{(ii)} and Propositions~\ref{prop:counterexample} and~\ref{prop:slow-quantization} show that, in several dimensions, the adapted quantization rate is governed by the stability of the diffusion flow, that is, by how strongly the flow separates nearby states.
Under strong stability, including the Lipschitz case, adapted quantization recovers the Brownian rate $(\log M)^{-1/2}$.
Under weaker stability, the rate deteriorates according to the Osgood flow, and Proposition~\ref{prop:counterexample} shows that this loss is genuine.
Without such quantitative stability, the rate can be slower than every negative power of $\log M$, even for smooth bounded coefficients.
This is in sharp contrast to classical functional quantization, which attains $(\log M)^{-1/2}$ for general It\^o processes under mild assumptions \citep{HeYuan2026}, including the examples of Proposition~\ref{prop:counterexample}, regardless of their stability.
The difference reflects what each problem has to approximate.
A classical quantizer only needs to approximate trajectories, whereas an adapted quantizer must also approximate, at every time, the conditional law of the future given the present.
A weakly stable flow can turn small differences in the present state into very different futures, which requires the quantizer to approximate more conditional laws.
\end{remark}

\section{Proofs}\label{sec:proofs}
\subsection{Proof strategy}
\label{sec:proof-strategy}

\paragraph{Theorem~\ref{thm:recursive-sampling}: Recursive sampling and upper bounds.}
To obtain a small relaxed adapted Wasserstein error, we need to construct a coupling with a small path error and a small delay in both causality directions. To this end, divide $[0,1]$ into intervals of length $h:=1/n$, with $t_k:=kh$, and {\it locally} construct the coupling between $X$ and uniform distribution over all sampled paths, one interval at a time. Lemma~\ref{lem:block-bicausal} is used to ensure $h$-bicausality.
This construction requires simultaneous control of the path error within each interval and the effect of endpoint discrepancies on subsequent evolution.

To be more specific, fix a sampling realization $\mathcal S$ and drive $X$ with a Brownian motion independent of the sampling.
Starting from $Y_0=x_0$, on each interval $[t_k,t_{k+1}]$ we generate a comparison diffusion from $Y_{t_k}$ using the same Brownian increment as $X$.
Each continuous path on the interval $[t_k,t_{k+1}]$ has a label recording both a normalized path codeword of accuracy $r$ from Lemma~\ref{lem:local-approximation} and an endpoint cell of accuracy $\delta$ from Lemma~\ref{lem:scalar-endpoints} or Lemma~\ref{lem:osgood-endpoints}. The comparison path (synchronously coupled with $X$) then determines a subset of samples through the maximal label coupling of Lemma~\ref{lem:maximal-label-coupling}. At every step, the matching uses only the updated path prefix and fresh interval randomness, so $Y$ has law $\mu^{n,L}$ and the recursion satisfies Lemma~\ref{lem:block-bicausal} ($h$-bicausality).
A successful match simultaneously controls the path quantization error within the interval and the endpoint discrepancy.
Lemma~\ref{lem:scalar-endpoints} controls the propagation of these discrepancies through scalar comparison principle, while Lemma~\ref{lem:osgood-endpoints} uses quantitative Osgood stability in higher dimensions.

Averaging over $\mathcal S$ recovers the original diffusion law, so $Y$ has the same unconditional moment bounds as $X$.
Let $G$ denote the event of successful path and endpoint matches on {\it every} interval, and let $\P$ include both sampling and coupling randomness.
Combining the local path estimates, the endpoint stability bound $\mathsf S_p(n\delta)$ defined in \eqref{eq:strategy-stability}, and the moment bounds on $G^c$, we obtain
\begin{equation}\label{eq:strategy-error}
 \mathcal E_{n,L,p}\leq\frac1n
 +C_p\left[\mathsf S_p(n\delta)+\frac r{\sqrt n}
 +\P(G^c)^{1/(2p)}\right].
\end{equation}
The four terms in \eqref{eq:strategy-error} account for the causality delay, endpoint error propagation, local path approximation, and matching failures. Choosing these parameters appropriately leads to Theorem~\ref{thm:recursive-sampling}.

\paragraph{Propositions~\ref{prop:counterexample} and~\ref{prop:slow-quantization}: Lower bounds.}
We derive both lower bounds from a common inequality.
Suppose that the filtration of $\mathbb X$ is constant on $[t_0,t_0+2\varepsilon_0]$, and let $R$ be a $1$-Lipschitz observation of the path.
Suppose further that there are laws $\mu_1,\ldots,\mu_J$ with pairwise $\cW_1$-distances at least $4r$, and disjoint events in $\cF^X_{t_0}$ with probabilities at least $w/J$, on which the conditional law of $R(X)$ given $\cF^X_{t_0}$ lies within $r$ of the corresponding $\mu_j$.
Then Corollary~\ref{cor:conditional-law-packing} gives
\begin{equation}\label{eq:strategy-lower}
 e^{\mathrm{ad}}_{M,1}(\mathbb X)\geq\min\{\varepsilon_0,wr/2\},\qquad 1\leq M\leq J/2.
\end{equation}
The bound in \eqref{eq:strategy-lower} is determined by the number $J$ of conditional laws, their separation $r$, and the probability $w$ with which the family is realized.
To obtain a single diffusion for all budgets, we combine countably many levels in the construction.
At level $k$, the diffusion produces $J_k$ laws with separation $r_k$ and reaches this family with probability of order $w_k$, yielding a lower bound for budgets near $M_k:=\lfloor J_k/2\rfloor$.
Monotonicity of $e^{\mathrm{ad}}_{M,1}$ in $M$ extends the bound to intermediate budgets.

Both constructions follow the same steps.
We first apply Lemma~\ref{lem:half-subset-packing} to obtain exponentially many separated uniform laws from finitely many separated points, and represent these laws by continuous samplers.
A diffusion selects a law before an information pause and samples it afterwards, so that its conditional law at the pause is the selected law.
Finally, we balance the parameters subject to the regularity constraint of the diffusion class.
The two examples differ in the observation and the regularity constraint on the samplers.

In the Osgood example of Proposition~\ref{prop:counterexample}, we observe the terminal value of $Z\in\R^q$ and construct the target laws on a grid of $2^{qk}$ points in $[a_k,2a_k]^q$.
The first Brownian stage selects a law, an independent increment after the pause supplies the sampling randomness, and the third block of coordinates $Z$ records and amplifies the sample.
The samplers that distinguish $e^{c2^{qk}}$ laws have large Lipschitz constants.
By composing these samplers with $g_m(z)=e^{-z^{-m}}$, we obtain globally Lipschitz seeds by choosing $a_k=A2^{-qk/m}$ with $A>0$ sufficiently small.
By Lemma~\ref{lem:osgood-amplifier}, applying $\Psi_1$ to each seed recovers the corresponding output through $\Psi_1(g_m(z))=z$.
Corollary~\ref{cor:osgood-vector-flow} extends this construction to all output coordinates and gives $\Phi(a)\asymp[\log(1/a)]^{-2/m}$, that is, $\beta=2/m$.
Taking $m=2/\beta$ and $q=\lceil m^2\rceil$, the grid spacing $r_k=a_k2^{-k}$ and layer masses of order $2^{-k}$ yield the lower exponent $1/m+2/q\leq\beta/2+\beta^2/2$ in dimension $d=q+2$.

In the slow rate example of Proposition~\ref{prop:slow-quantization}, we observe the entire path of a scalar coordinate $Z$ and use the $2^{n_k}$ binary paths obtained by selecting bumps on $n_k$ disjoint time intervals as the separated points.
One coordinate selects a cell before the pause and supplies the sampling randomness afterwards, while the second coordinate records the sampled path by integrating its time derivative.
The samplers require only a bounded time derivative, which allows us to choose the bump height $r_k=1/n_k$.
We place the levels in disjoint cells tending to spatial infinity to preserve smoothness of the coefficients, at the cost of the Gaussian probability $e^{-(k+1)^2/2}$ of reaching level $k$.
With $n_k=k^3$, the growth $\log M_k\asymp2^{k^3}$ compensates for this probability decay, and the lower bound decays more slowly than every negative power of $\log M$.

\subsection{Technical lemmas for the upper bounds}
\label{sec:stability}

We collect the lemmas used in the proof of Theorem~\ref{thm:recursive-sampling}, beginning with a general principle for constructing relaxed bicausal couplings from Brownian increments. Throughout this subsection,
fix $h:=1/n$, $t_k:=kh$, and let $\mathcal G$ be the usual augmentation of the Brownian filtration.
The next lemma allows each segment to depend on the full observed prefix and independent auxiliary randomness.
In the recursive sampling construction, we apply it with the sampling realization $\mathcal S$ fixed.

\begin{lemma}\label{lem:block-bicausal}
Suppose Assumption~\ref{ass:coefficients} holds.
Let $U_0,\ldots,U_{n-1}$ be independent uniform random variables on $[0,1]$, independent of the entire Brownian path $B$.
For $0\leq k<n$, let
\begin{equation}
 F_k:C([0,t_k];\R^d)\times\cC_m\times[0,1]\longrightarrow\cC_d
\end{equation}
be a Borel map satisfying $F_k(g,w,u)(0)=g(t_k)$.
Define $Y_0:=x_0$ and
\begin{equation}
 (Y_{t_k+uh})_{0\leq u\leq1}
 :=F_k\bigl(Y|_{[0,t_k]},(B_{t_k+uh}-B_{t_k})_{0\leq u\leq1},U_k\bigr).
\end{equation}
Then $\Lc(X,Y)$ is an $h$-bicausal coupling of $\mathbb X$ and $\Can(\Lc(Y))$.
\end{lemma}

\begin{proof}
Let $\mathcal H$ be the usual augmentation of the filtration
\begin{equation}
 \bigl(\sigma(B_v:0\leq v\leq t)
 \vee\sigma(U_j:t_{j+1}\leq t)\bigr)_{0\leq t\leq1}.
\end{equation}
Thus $U_j$ is revealed at $t_{j+1}$, and $B$ remains Brownian with respect to $\mathcal H$.
We complete all $\sigma$-fields below modulo the joint law's null sets.
Strong well-posedness and independence of future Brownian increments give, for every bounded Borel function $f$ on path space,
\begin{equation}\label{eq:immersion}
 \E[f(X)\mid\mathcal H_t]=\E[f(X)\mid\cF_t^X],\qquad 0\leq t\leq1.
\end{equation}

Let $\mathcal H_t^Y:=\sigma(Y_v:0\leq v\leq t)$, completed as above. Note that in our construction, $\{\mathcal H^Y_t\}_{0\leq t\leq 1}$ is not necessarily right-continuous, and $\{\mathcal F^Y_t\}_{0\leq t\leq 1}$, the canonical filtration, is its right-continuous augmentation. At each grid time $t_j$, the recursion gives $\mathcal H_{t_j}^Y\subseteq\mathcal H_{t_j}$.
The future path of $Y$ is a Borel function of $Y|_{[0,t_j]}$, the Brownian increments after $t_j$, and $U_j,\ldots,U_{n-1}$.
These future noises are independent of $\mathcal H_{t_j}$, so
\begin{equation}\label{eq:block-prefix-conditional}
 \E[f(Y)\mid\mathcal H_{t_j}]=\E[f(Y)\mid\mathcal H_{t_j}^Y].
\end{equation}

Both causality conditions are immediate when $t+h\geq1$.
We therefore fix $t\in[t_k,t_{k+1})$ with $t+h<1$.
For every $u<t_{k+1}$, the prefix $Y|_{[0,u]}$ is $\mathcal H_{t_{k+1}}$-measurable.
Since $t<t_{k+1}$, taking right limits through times $u\in(t,t_{k+1})$ therefore gives
\begin{equation}
 \cF_t^Y\subseteq\mathcal H_{t_{k+1}}\subseteq\mathcal H_{t+h},
 \qquad \cF_t^X\subseteq\mathcal H_{t_{k+1}},
 \qquad \mathcal H_{t_{k+1}}^Y\subseteq\cF_{t+h}^Y.
\end{equation}
These inclusions also hold when $t$ is a grid time.
For bounded, $\cF_t^Y$-measurable $V$, condition first on $\mathcal H_{t+h}$ and then on $\cF_{t+h}^X$.
By \eqref{eq:immersion},
\begin{equation}
 \E[Vf(X)\mid\cF_{t+h}^X]
 =\E[V\mid\cF_{t+h}^X]\E[f(X)\mid\cF_{t+h}^X].
\end{equation}
This proves causality from $X$ to $Y$.

For reverse causality, let $U$ be bounded and $\cF_t^X$-measurable, hence $\mathcal H_{t_{k+1}}$-measurable.
For every bounded Borel $f$, \eqref{eq:block-prefix-conditional} yields
\begin{equation}
 \E[Uf(Y)]
 =\E\bigl[U\E[f(Y)\mid\mathcal H_{t_{k+1}}^Y]\bigr]
 =\E\bigl[\E[U\mid\mathcal H_{t_{k+1}}^Y]f(Y)\bigr].
\end{equation}
Consequently, using $\mathcal H_{t_{k+1}}^Y\subseteq\cF_{t+h}^Y\subseteq\cF_1^Y$ and the tower property,
\begin{equation}
 \E[U\mid\cF_1^Y]
 =\E[U\mid\mathcal H_{t_{k+1}}^Y]
 =\E[U\mid\cF_{t+h}^Y].
\end{equation}
This is reverse causality.
\end{proof}

The next lemma combines the standard maximal coupling on a finite set \citep[Proposition~4.7]{LevinPeres2017} with an empirical frequency estimate. It will be used to construct the $h$-bicausal coupling between $X$ and $Y$, and part (ii) also provides a bound on mismatch probability, which is important in the estimate of failure probability in the main proof.

\begin{lemma}\label{lem:maximal-label-coupling}
Let $E$ be a standard Borel space and $\ell:E\to\{1,\ldots,A\}$ a Borel label map, where $A\geq1$.
\begin{enumerate}[label=\textnormal{(\roman*)}]
\item Let $P$ be a probability measure on $E$ and $s_1,\ldots,s_N\in E$ a fixed list, $N\geq1$, counting repetitions separately.
Set
\begin{equation}
 \lambda_a:=P(\ell=a),\qquad
 \widehat\lambda_a:=\frac1N\sum_{j=1}^N\1_{\{\ell(s_j)=a\}},\qquad 1\leq a\leq A.
\end{equation}
Given $C\sim P$ and an independent uniform random variable $U$ on $[0,1]$, one can select $J\in\{1,\ldots,N\}$ by a Borel function of $(\ell(C),\ell(s_1),\ldots,\ell(s_N),\lambda,U)$ such that
\begin{equation}\label{eq:maximal-label-matching}
 \P(J=j)=\frac1N,\qquad
 \P\bigl(\ell(C)\ne\ell(s_J)\bigr)
 =\TV(\lambda,\widehat\lambda)
 :=\frac12\sum_{a=1}^A|\lambda_a-\widehat\lambda_a|.
\end{equation}
\item Suppose that, conditionally on a $\sigma$-field $\mathcal D$, $S_1,\ldots,S_N$ are i.i.d.\ with law $P$, where $N\geq1$ and the probability kernel $P$ on $E$ are $\mathcal D$-measurable.
Let $C$ have conditional law $P$ independently of the sample list given $\mathcal D$, and select $J$ by \textnormal{(i)} using a uniform $U$ on $[0,1]$ independent of $\mathcal D\vee\sigma(C,S_1,\ldots,S_N)$.
Then
\begin{equation}\label{eq:empirical-label-tv}
 \E\!\left[\P\bigl(\ell(C)\ne\ell(S_J)\mid\mathcal D,S_1,\ldots,S_N\bigr)\,\middle|\,\mathcal D\right]
 \leq\frac12\sqrt{\frac{A-1}{N}}\qquad\text{almost surely}.
\end{equation}
\end{enumerate}
\end{lemma}

\begin{proof}
\textnormal{(i)} We first couple the labels, then distribute the selected label's mass uniformly among its sample indices.
Set $m_a:=\min\{\lambda_a,\widehat\lambda_a\}$ and $D:=1-\sum_a m_a=\TV(\lambda,\widehat\lambda)$.
For $D>0$, define
\begin{equation}
 \Gamma_{ab}:=m_a\1_{\{a=b\}}
 +\frac{(\lambda_a-m_a)(\widehat\lambda_b-m_b)}{D},
 \qquad 1\leq a,b\leq A.
\end{equation}
For $D=0$, set $\Gamma_{ab}:=\lambda_a\1_{\{a=b\}}$.
The row and column sums of $\Gamma$ are $\lambda$ and $\widehat\lambda$.
The residual masses $\lambda_a-m_a$ and $\widehat\lambda_a-m_a$ cannot both be positive, so $\Gamma_{aa}=m_a$ and $\sum_{a\ne b}\Gamma_{ab}=D$.

Write $N_a:=\#\{j:\ell(s_j)=a\}$.
Given $C=c$ with $\lambda_{\ell(c)}>0$, select an index according to
\begin{equation}
 K(c,j):=\frac{\Gamma_{\ell(c),\ell(s_j)}}{\lambda_{\ell(c)}N_{\ell(s_j)}},
 \qquad 1\leq j\leq N.
\end{equation}
Each $N_{\ell(s_j)}$ is positive.
Grouping indices by label and using the row sums of $\Gamma$ shows that $K(c,\cdot)$ is a probability vector.
On the $P$-null set where $\lambda_{\ell(c)}=0$, set $K(c,j):=1/N$.
The entries of $K$ are Borel functions of the labels and $\lambda$, so the inverse cumulative distribution rule using $U$ gives the required Borel selection.
For every $a$ and $j$, this construction gives
\begin{equation}
 \P\bigl(\ell(C)=a,J=j\bigr)
 =\frac{\Gamma_{a,\ell(s_j)}}{N_{\ell(s_j)}}.
\end{equation}
Summing over $a$ yields $\P(J=j)=\widehat\lambda_{\ell(s_j)}/N_{\ell(s_j)}=1/N$, while summing over indices of each label shows that $(\ell(C),\ell(s_J))$ has law $\Gamma$.
Its off-diagonal mass is $D$, which proves \eqref{eq:maximal-label-matching}.

\textnormal{(ii)} Write $\lambda_a:=P(\ell=a)$ and let $\widehat\lambda_N$ be the empirical label law of $S_1,\ldots,S_N$.
Conditionally on $\mathcal D$ and this list, part~\textnormal{(i)} gives
\begin{equation}
 \P\bigl(\ell(C)\ne\ell(S_J)\mid\mathcal D,S_1,\ldots,S_N\bigr)
 =\TV(\lambda,\widehat\lambda_N).
\end{equation}
Conditionally on $\mathcal D$, the sample labels are independent with law $\lambda$.
Hence, for every $a$,
\begin{equation}
 \E\bigl[\widehat\lambda_{N,a}\mid\mathcal D\bigr]=\lambda_a,\qquad
 \E\bigl[(\widehat\lambda_{N,a}-\lambda_a)^2\mid\mathcal D\bigr]
 =\frac{\lambda_a(1-\lambda_a)}N.
\end{equation}
Conditional Cauchy--Schwarz therefore gives
\begin{equation}
 \E\bigl[\TV(\lambda,\widehat\lambda_N)\mid\mathcal D\bigr]
 \leq\frac1{2\sqrt N}\sum_{a=1}^A\sqrt{\lambda_a(1-\lambda_a)}
 \leq\frac1{2\sqrt N}\sqrt{A\left(1-\sum_{a=1}^A\lambda_a^2\right)}
 \leq\frac12\sqrt{\frac{A-1}{N}},
\end{equation}
where the last inequality uses $\sum_a\lambda_a^2\geq1/A$.
This proves \eqref{eq:empirical-label-tv}.
\end{proof}

To construct the $h$-bicausal coupling between $X$ and $Y$ through labels, we use the deterministic codebook of \citet[Proposition~4.2]{HeYuan2026}. To apply the codebook estimate under linear growth, we normalize each short segment by its random initial size and use a stopped version with bounded coefficients.
The following lemma gives such a codebook along with a conditional bound on the coding error.
Here $\mathcal G$ may be any filtration satisfying the usual conditions for which $B$ is Brownian.

\begin{lemma}\label{lem:local-approximation}
Fix $d\in\N$.
Let $K\geq1$ be arbitrary, and let $Z$ be any $\R^d$-valued It\^o process of the form
\begin{equation}
 Z_t=Z_0+\int_0^t a_s\,\mathrm ds+\int_0^t H_s\,\mathrm dB_s,\qquad 0\leq t\leq h,
\end{equation}
where $Z_0\in\R^d$ is $\mathcal G_0$-measurable, $B$ is a $\mathcal G$-Brownian motion, and $a,H$ are progressively measurable with respect to $\mathcal G$ and satisfy
\begin{equation}
 |a_t|+\|H_t\|_{\mathrm{HS}}\leq K(1+|Z_t|),\qquad 0\leq t\leq h.
\end{equation}
Then, for all sufficiently small $h,r>0$, there exist a deterministic nonempty finite set $\cA_r\subset\cC_d$ and a Borel map $\alpha_{h,r}:\cC_d\to\cA_r$ such that
\begin{equation}\label{eq:forest-size}
 \log|\cA_r|\leq C_dK^2r^{-2},
\end{equation}
and the normalized path $V_u:=(Z_{hu}-Z_0)/[\sqrt h(1+|Z_0|)]$, $0\leq u\leq1$, satisfies
\begin{equation}\label{eq:local-approximation}
 \Pp\left(\|V-\alpha_{h,r}(V)\|_\infty>C_dr\,\middle|\,\mathcal G_0\right)\leq e^{-c_dK^2/r^2}+Ce^{-c/h}
 \qquad\text{almost surely}.
\end{equation}
The set $\cA_r$ depends only on $d,K,r$, and the map $\alpha_{h,r}$ depends only on $d,K,h,r$.
The constants $c_d,C_d>0$ depend only on $d$, while $c,C>0$ and the thresholds depend only on $d,K$.
\end{lemma}

\begin{proof}
Set $W_u:=B_{hu}/\sqrt h$ and
\begin{equation}
 \tau:=\inf\{u\in[0,1]:|V_u|\geq h^{-1/2}\}\wedge1,
 \qquad \widetilde V_u:=V_{u\wedge\tau}.
\end{equation}
Then $W$ is Brownian for the filtration $(\mathcal G_{hu})_{0\leq u\leq1}$.
For $u<\tau$, we have $1+|Z_{hu}|\leq2(1+|Z_0|)$, and hence
\begin{equation}
 \frac{\sqrt h\,|a_{hu}|}{1+|Z_0|}\leq2K\sqrt h,
 \qquad \frac{\|H_{hu}\|_{\mathrm{HS}}}{1+|Z_0|}\leq2K.
\end{equation}
Both integrands of $\widetilde V$ vanish after stopping.
For $h\leq1$, its drift and diffusion integrands are therefore bounded by $2K$.
Apply the deterministic codebook construction of \citet[Proposition~4.2]{HeYuan2026} with $T=1$, $G=2K$, and approximation scale $r$, fixing the initial lattice point at zero since $\widetilde V_0=0$.
The resulting nonempty codebook $\cA_r$ depends only on $d,K,r$ and satisfies \eqref{eq:forest-size}.

Fix an enumeration of $\cA_r$ and let $\alpha_{h,r}(w)$ be its first element minimizing $\|w-a\|_\infty$ over $a\in\cA_r$.
The map $\alpha_{h,r}$ is Borel because the codebook is finite.
For every $A\in\mathcal G_0$ with $\Pp(A)>0$, $W$ remains Brownian under $\Pp(\,\cdot\mid A)$ because its increments are independent of the past, including $A$.
The same deterministic codebook and coefficient bounds apply under each such measure.
Since $A$ is arbitrary, the cited error-probability estimate gives
\begin{equation}
 \Pp\left(\|\widetilde V-\alpha_{h,r}(\widetilde V)\|_\infty>C_dr\,\middle|\,\mathcal G_0\right)
 \leq e^{-c_dK^2/r^2}\qquad\text{almost surely}.
\end{equation}
It remains to bound the probability of reaching the stopping threshold.
The drift of $\widetilde V$ has total variation at most $2K\sqrt h$, and each coordinate of its martingale part $M$ has quadratic variation at most $4K^2$.
On $\{\|V\|_\infty\geq h^{-1/2}\}$, continuity gives $\|\widetilde V\|_\infty=h^{-1/2}$.
For $h\leq(4K)^{-1}$, the drift contributes at most half this threshold.
Thus the conditional martingale maximal estimate in \citet[Lemma~4.1]{HeYuan2026}, applied to each coordinate, gives
\begin{equation}
 \Pp\left(\|V\|_\infty\geq h^{-1/2}\,\middle|\,\mathcal G_0\right)
 \leq\sum_{i=1}^d\Pp\left(\sup_{u\leq1}|M_u^i|\geq\frac{1}{2\sqrt{dh}}\,\middle|\,\mathcal G_0\right)\leq2d\,e^{-1/(32dK^2h)}.
\end{equation}

On $\{\|V\|_\infty<h^{-1/2}\}$, we have $\widetilde V=V$.
Thus
\begin{equation}
 \{\|V-\alpha_{h,r}(V)\|_\infty>C_dr\}
 \subseteq\{\|V\|_\infty\geq h^{-1/2}\}
 \cup\{\|\widetilde V-\alpha_{h,r}(\widetilde V)\|_\infty>C_dr\},
\end{equation}
and the two conditional probability bounds give \eqref{eq:local-approximation}.
\end{proof}

For the comparison segment $C^{(k)}$ on $[t_k,t_{k+1}]$, the proof of Theorem~\ref{thm:recursive-sampling} applies Lemma~\ref{lem:local-approximation} to $Z_u:=C^{(k)}_{t_k+u}$, $0\leq u\leq h$, driven by $B_{t_k+u}-B_{t_k}$ in the shifted filtration $(\mathcal H_{t_k+u})_{0\leq u\leq h}$, where $\mathcal H$ is as in the proof of Lemma~\ref{lem:block-bicausal}.
Its initial state $Z_0=Y_{t_k}$ is $\mathcal H_{t_k}$-measurable, so the coding-error bound holds conditionally on $\mathcal H_{t_k}$.

We now establish the stability estimates for accumulated endpoint perturbations. Fix $n\in\N$ and set $t_k:=k/n$, $k=0,\ldots,n$.
Let $(\mathcal H_t)_{0\leq t\leq1}$ be a usual filtration in which $B$ remains Brownian.
Let $Z$ be an $\mathcal H$-adapted c\`adl\`ag process with $Z_0=x_0$ which, between grid times, solves the original SDE driven by the same $B$ as $X$, and has jumps only at grid times. In other words, $Z$ satisfies
\begin{equation}\label{eq:quantized-recursion}
 Z_t=x_0+\int_0^t b(s,Z_{s-})\,\mathrm ds
 +\int_0^t\sigma(s,Z_{s-})\,\mathrm dB_s
 +\sum_{j:\,t_j\leq t}\Delta Z_{t_j},\qquad 0\leq t\leq1.
\end{equation}

\begin{lemma}\label{lem:scalar-endpoints}
Suppose $d=1$.
Set $v(t,x):=\E\int_t^1X_s^{t,x}\,\mathrm ds$ and $f_t(x):=x+v(t,x)$.
For each $t$, the function $v(t,\cdot)$ is continuous and nondecreasing, with $|v(t,x)|\leq C(1+|x|)$ uniformly in $t$.
Moreover, $f_t$ is a continuous strictly increasing bijection of $\R$, and $f_t^{-1}$ is $1$-Lipschitz.
If $\delta>0$ satisfies $n\delta\leq1$ and the process $Z$ in \eqref{eq:quantized-recursion} satisfies
\begin{equation}\label{eq:scalar-adapted-jumps}
 |f_{t_j}(Z_{t_j})-f_{t_j}(Z_{t_j-})|\leq\delta,\qquad j=1,\ldots,n,
\end{equation}
then, for every $p\geq1$,
\begin{equation}\label{eq:scalar-endpoint-bound}
 \left(\E\sup_{t\leq1}|X_t-Z_t|^p\right)^{1/p}\leq C_p(n\delta)^{1/[5(2p-1)]}.
\end{equation}
\end{lemma}

\begin{proof}
For each fixed starting time $t$ and pair $x\leq y$, couple $X^{t,x}$ and $X^{t,y}$ with the same Brownian motion.
Continuity and pathwise uniqueness prevent these scalar solutions from crossing, so $X_s^{t,x}\leq X_s^{t,y}$ for every $s\in[t,1]$ almost surely.
Order preservation and uniform linear growth give $v(t,x)\leq v(t,y)$ for $x\leq y$ and $|v(t,x)|\leq C(1+|x|)$, uniformly in $t$.
For fixed $t$ and $x_j\to x$, the moment and increment bounds give tightness of $X^{t,x_j}$ on $C([t,1];\R)$.
Continuity of the coefficients identifies each subsequential limit as a weak solution from $(t,x)$, and uniqueness in law therefore gives $X^{t,x_j}\Rightarrow X^{t,x}$.
Uniform moment bounds then imply continuity of $v(t,\cdot)$.
Hence $f_t$ is continuous and strictly increasing, with $f_t(x)-f_t(y)\geq x-y$ for $x\geq y$.
It tends to the corresponding infinity at both ends of $\R$, so it is onto and its inverse is $1$-Lipschitz.
In particular, \eqref{eq:scalar-adapted-jumps} implies $|\Delta Z_{t_j}|\leq\delta$.
The total jump variation of $Z$ is therefore at most $n\delta\leq1$.
Linear growth, the Burkholder--Davis--Gundy inequality, and Gronwall's inequality give
\begin{equation}\label{eq:endpoint-jump-moments}
 \E\sup_{t\leq1}|Z_t|^q+\E\sup_{t\leq1}|X_t-Z_t|^q\leq C_q,\qquad q\geq1,
\end{equation}
uniformly over $n$, $\delta$, and all jumps satisfying \eqref{eq:scalar-adapted-jumps}.
These bounds also control the left limits at grid times.

To compare $Z$ with $X$, define the increasing maps
\begin{equation}
 R_{t,\delta}^{\pm}(x):=f_t^{-1}(f_t(x)\pm\delta).
\end{equation}
Since $f_t=x+v(t,x)$ with both terms increasing,
\begin{align}
 0&\leq R_{t,\delta}^+(x)-x\leq\delta,
 &0&\leq v(t,R_{t,\delta}^+(x))-v(t,x)\leq\delta,\notag\\
 0&\leq x-R_{t,\delta}^-(x)\leq\delta,
 &0&\leq v(t,x)-v(t,R_{t,\delta}^-(x))\leq\delta.
\end{align}
Let $Z^+$ and $Z^-$ start at $x_0$, evolve between grid times with the same $B$ as $X$ and $Z$, and reset at $t_j$ by $R_{t_j,\delta}^+$ and $R_{t_j,\delta}^-$, respectively.
Condition \eqref{eq:scalar-adapted-jumps} gives
$R_{t_j,\delta}^-(Z_{t_j-})\leq Z_{t_j}\leq R_{t_j,\delta}^+(Z_{t_j-})$.
At each restart the states are $\mathcal H_{t_j}$-measurable, and the future increments of $B$ are independent of $\mathcal H_{t_j}$.
Thus pathwise uniqueness and scalar comparison apply from these random initial states.
Comparison between grid times, followed by monotonicity of $R^\pm$, yields
\begin{equation}
 Z_t^-\leq\min\{X_t,Z_t\}\leq\max\{X_t,Z_t\}\leq Z_t^+,\qquad 0\leq t\leq1.
\end{equation}
Both envelopes satisfy \eqref{eq:endpoint-jump-moments}, since each of their jumps has size at most $\delta$.

For $W=Z^+$ or $Z^-$, condition on $\mathcal H_{t_{j-1}}$ and apply the Markov property to the diffusion started at $W_{t_{j-1}}$.
The tower property then gives
\begin{equation}
 \E\int_{t_{j-1}}^{t_j}W_s\,\mathrm ds
 =\E\bigl[v(t_{j-1},W_{t_{j-1}})-v(t_j,W_{t_j-})\bigr].
\end{equation}
Summing over $j$, using $v(1,\cdot)=0$ and $\E\int_0^1X_s\,\mathrm ds=v(0,x_0)$, we obtain
\begin{align}
 \E\int_0^1(Z_s^+-X_s)\,\mathrm ds
 &=\sum_{j=1}^n\E\bigl[v(t_j,Z_{t_j}^+)-v(t_j,Z_{t_j-}^+)\bigr]\leq n\delta,\notag\\
 \E\int_0^1(X_s-Z_s^-)\,\mathrm ds
 &=\sum_{j=1}^n\E\bigl[v(t_j,Z_{t_j-}^-)-v(t_j,Z_{t_j}^-)\bigr]\leq n\delta.
 \label{eq:integrated-scalar-error}
\end{align}

To pass from these integral bounds to a uniform bound, take either nonnegative process $D=Z^+-X$ or $D=X-Z^-$ and let
\begin{equation}
 A_t:=D_t-\sum_{j:\,t_j\leq t}\Delta D_{t_j}
\end{equation}
be its continuous semimartingale part.
The moment bounds for $X$ and $Z^\pm$, linear growth, and the Burkholder--Davis--Gundy inequality give $\|A_t-A_s\|_{L^q}\leq C_q|t-s|^{1/2}$ for $s,t\in[0,1]$ and $q\geq2$.
Taking $q=8$, we have $\E|A_t-A_s|^8\leq C|t-s|^4$.
The dyadic chaining estimate in the proof of Kolmogorov's continuity theorem \citep[Chapter~I, Section~2, Theorem~(2.1)]{RevuzYor1999} therefore gives, for $0<\eta\leq1$,
\begin{equation}
 \E\sup_{|t-s|\leq\eta}|A_t-A_s|
 \leq C\eta^{1/4}.
\end{equation}

The total jump variation of $D$ is at most $n\delta$, so, for all $s,t\in[0,1]$,
\begin{equation}
 |D_t-D_s|\leq|A_t-A_s|+n\delta.
\end{equation}
For each $t$, choose $I_t\subset[0,1]$ of length $\eta$ containing $t$.
Averaging over $s\in I_t$ gives
\begin{equation}
 D_t\leq\eta^{-1}\int_{I_t}D_s\,\mathrm ds
       +\sup_{|u-v|\leq\eta}|A_u-A_v|+n\delta.
\end{equation}
Since $D\geq0$, enlarge the integral to $[0,1]$ and take the supremum over $t$.
By \eqref{eq:integrated-scalar-error},
\begin{equation}
 \E\sup_{t\leq1}D_t\leq\frac{n\delta}{\eta}+C\eta^{1/4}+n\delta.
\end{equation}
Choose $\eta:=(n\delta)^{4/5}$.
As $n\delta\leq1$, each envelope difference satisfies $\E\sup_tD_t\leq C(n\delta)^{1/5}$.
The envelope inequality gives $|X_t-Z_t|\leq(Z_t^+-X_t)+(X_t-Z_t^-)$, hence
\begin{equation}
 \E\sup_{t\leq1}|X_t-Z_t|\leq C(n\delta)^{1/5}.
\end{equation}
For $p>1$, H\"older's inequality and the $2p$-moment bound in \eqref{eq:endpoint-jump-moments} give
\begin{align}
 \left(\E\sup_t|X_t-Z_t|^p\right)^{1/p}
 &\leq\left(\E\sup_t|X_t-Z_t|\right)^{1/(2p-1)}
 \times\left(\E\sup_t|X_t-Z_t|^{2p}\right)^{(p-1)/[p(2p-1)]}\notag\\
 &\leq C_p(n\delta)^{1/[5(2p-1)]}.
\end{align}
The first-moment estimate gives the same assertion for $p=1$, proving \eqref{eq:scalar-endpoint-bound}.
\end{proof}

In several dimensions, the Osgood condition controls the error caused by bounded adapted jumps.

\begin{lemma}\label{lem:osgood-endpoints}
Suppose Assumptions~\ref{ass:coefficients} and~\ref{ass:osgood} hold.
Fix $c_d>0$ and let $Z$ be the process in \eqref{eq:quantized-recursion}.
If $\delta>0$ satisfies $n\delta\leq1$ and
\begin{equation}
 |\Delta Z_{t_j}|\leq c_d\delta\quad\text{almost surely},\qquad j=1,\ldots,n,
\end{equation}
then, for every $p\geq1$,
\begin{equation}
 \left(\E\sup_{t\leq1}|X_t-Z_t|^p\right)^{1/p}
 \leq C_p\Phi(C_pn\delta)^{1/(4p-2)}.
\end{equation}
The constant $C_p$ may depend on $p$, $c_d$, and the data in the assumptions, but is independent of $n$, $\delta$, and the choice of adapted jumps.
Here $\Phi$ is the Osgood flow in \eqref{eq:comparison-flow}.
\end{lemma}

\begin{proof}
The total jump variation of $Z$ is at most $c_dn\delta\leq c_d$.
The moment argument for \eqref{eq:endpoint-jump-moments} in the proof of Lemma~\ref{lem:scalar-endpoints} uses only this variation bound, linear growth, the Burkholder--Davis--Gundy inequality, and Gronwall's inequality, so it applies in any dimension.
Thus \eqref{eq:endpoint-jump-moments} holds for the present $Z$ and also controls its left limits, uniformly over all admissible jumps.

To estimate the error, write $D_t:=X_t-Z_t$ and, for $v>0$, define the $\mathcal H$-stopping time
\begin{equation}
 \tau_v:=\inf\{t\in[0,1]:|D_t|>v\}\wedge1,\qquad \inf\varnothing:=+\infty.
\end{equation}
Between grid times, \eqref{eq:joint-osgood} bounds the drift of $|D|^2$ by $\ell(s)|D_s|^2+k(s)\rho(|D_s|^2)$.
At each grid time,
\begin{equation}
 |D_{t_k}|^2-|D_{t_k-}|^2
 =-2\langle D_{t_k-},\Delta Z_{t_k}\rangle+|\Delta Z_{t_k}|^2
 \leq2c_d\delta|D_{t_k-}|+c_d^2\delta^2.
\end{equation}
Consequently, stopped It\^o's formula gives
\begin{equation}
 |D_{t\wedge\tau_v}|^2
 \leq\int_0^{t\wedge\tau_v}\bigl[\ell(s)|D_s|^2+k(s)\rho(|D_s|^2)\bigr]\,\mathrm ds+M_{t\wedge\tau_v}
 +\sum_{j:\,t_j\leq t\wedge\tau_v}\bigl(2c_d\delta|D_{t_j-}|+c_d^2\delta^2\bigr),
\end{equation}
where
\begin{equation}
 M_t:=2\int_0^t\langle D_{s-},[\sigma(s,X_s)-\sigma(s,Z_{s-})]\,\mathrm dB_s\rangle.
\end{equation}
The fourth-moment bounds and linear growth imply $\E\langle M\rangle_1<\infty$, so $\E M_{t\wedge\tau_v}=0$.
Enlarging the nonnegative jump sum to all $n$ grid times bounds its expectation by
\begin{equation}
 2c_d\delta\sum_{j=1}^n\E|D_{t_j-}|+c_d^2n\delta^2\leq Cn\delta.
\end{equation}
Taking expectations and using the nonnegativity of $\ell,k,\rho$ and Jensen's inequality for the concave function $\rho$, we obtain
\begin{equation}
 \E|D_{t\wedge\tau_v}|^2
 \leq Cn\delta+\int_0^t\bigl[\ell(s)\E|D_{s\wedge\tau_v}|^2
 +k(s)\rho(\E|D_{s\wedge\tau_v}|^2)\bigr]\,\mathrm ds.
\end{equation}
Since $\rho$ is increasing and locally Lipschitz on $(0,\infty)$, the integral form of the scalar ODE comparison theorem applied to \eqref{eq:comparison-flow} gives
\begin{equation}
 \E|D_{t\wedge\tau_v}|^2\leq z(t;Cn\delta),\qquad 0\leq t\leq1.
\end{equation}
In particular, $\E|D_{\tau_v}|^2\leq\Phi(Cn\delta)$.

To pass to a uniform bound, observe that on $\{\sup_t|D_t|>v\}$ we have $|D_{\tau_v}|\geq v$. Thus
\begin{equation}
 \Pp(\sup_t|D_t|>v)\leq\min\{1,\Phi(Cn\delta)/v^2\},
\end{equation}
and integration over $v>0$ yields
\begin{equation}
 \E\sup_t|D_t|\leq2\Phi(Cn\delta)^{1/2}.
\end{equation}
For $p>1$, interpolate this first-moment estimate with $\E\sup_t|D_t|^{2p}\leq C_p$ as in the proof of Lemma~\ref{lem:scalar-endpoints}:
\begin{equation}
 \left(\E\sup_t|D_t|^p\right)^{1/p}
 \leq C_p\left(\E\sup_t|D_t|\right)^{1/(2p-1)}
 \leq C_p\Phi(Cn\delta)^{1/(4p-2)}.
\end{equation}
The first-moment estimate already gives the assertion for $p=1$.
\end{proof}

\subsection{Technical lemmas for the lower bounds}
\label{sec:lower-lemmas}

We divide the results used for the lower bounds into two groups.
The first group converts separated conditional laws into lower bounds.
On an interval where the filtration of $\mathbb X$ is constant, Lemma~\ref{lem:conditional-law-lower} compares $\mathbb X$ with a finitely supported law through the conditional laws of a Lipschitz observation of the path.
We use a terminal value as the observation in the proof of Proposition~\ref{prop:counterexample} and the path of one coordinate in the proof of Proposition~\ref{prop:slow-quantization}.
Corollary~\ref{cor:conditional-law-packing} gives a lower bound for the resulting approximation error when the conditional laws approximate well-separated target laws, and Lemma~\ref{lem:half-subset-packing} constructs exponentially many target laws from finitely many separated points.
The second group consists of Lemma~\ref{lem:osgood-amplifier} and Corollary~\ref{cor:osgood-vector-flow}, which provide the amplification mechanism for the Osgood example.

\begin{lemma}\label{lem:conditional-law-lower}
Let $p\geq1$ and let $\mathbb X$ be a canonically filtered process on $\cC_d$ with $\E\|X\|_\infty^p<\infty$.
Suppose there exist constants $t_0\in[0,1)$ and $\varepsilon_0>0$, with $t_0+2\varepsilon_0\leq1$, such that
\[
 \cF_t^X=\cF_{t_0}^X\qquad\text{for every }t\in[t_0,t_0+2\varepsilon_0].
\]
Fix such a pair $(t_0,\varepsilon_0)$.
Let $(E,d_E)$ be Polish and let $R:\cC_d\to E$ be $1$-Lipschitz.
Set $\mathcal K:=\Lc(R(X)\mid\cF_{t_0}^X)$ and
\begin{equation}
 Q_{M,p}:=\inf_{\lambda_1,\ldots,\lambda_M\in\cP_p(E)}
 \left(\E\min_{j\leq M}\cW_p(\mathcal K,\lambda_j)^p\right)^{1/p},
\end{equation}
where $\cW_p$ uses the distance $d_E$.
Then, for every $M\in\N$,
\begin{equation}\label{eq:conditional-law-lower}
 e^{\mathrm{ad}}_{M,p}(\mathbb X)\geq\min\{\varepsilon_0,Q_{M,p}\}.
\end{equation}
\end{lemma}

\begin{proof}
Fix an $\varepsilon$-bicausal coupling of $X,Y$ with $0\leq\varepsilon\leq\varepsilon_0$ and $|\supp\Lc(Y)|\leq M$.
Put $\mathcal H:=\cF_{t_0}^X\vee\cF_{t_0+\varepsilon}^Y$.
Since $\cF_{t_0+2\varepsilon}^X=\cF_{t_0}^X$, the two causality conditions at times $t_0+\varepsilon$ and $t_0$, respectively, give
\begin{equation}
 \Lc(R(X)\mid\mathcal H)=\Lc(R(X)\mid\cF_{t_0}^X),
 \qquad 
 \Lc(R(Y)\mid\mathcal H)=\Lc(R(Y)\mid\cF_{t_0+\varepsilon}^Y).
\end{equation}
Hence the joint law of $(R(X),R(Y))$ conditional on $\mathcal H$ couples
$\mathcal K$ and $\nu:=\Lc(R(Y)\mid\cF_{t_0+\varepsilon}^Y)$.
The Lipschitz property of $R$ yields
\begin{equation}
 \E[\|X-Y\|_\infty^p\mid\mathcal H]\geq \cW_p(\mathcal K,\nu)^p.
\end{equation}
On the support of $\Lc(Y)$, its canonical $\sigma$-field at time $t_0+\varepsilon$ has at most $M$ atoms, modulo null sets.
Thus $\nu$ takes at most $M$ values $\lambda_1,\ldots,\lambda_M\in\cP_p(E)$, and
\begin{equation}
 \E\|X-Y\|_\infty^p
 \geq\E\min_{j\leq M}\cW_p(\mathcal K,\lambda_j)^p
 \geq Q_{M,p}^p.
\end{equation}
Every coupling with $\varepsilon>\varepsilon_0$ has relaxed cost at least $\varepsilon_0$.
Taking the infima in \eqref{eq:aw-definition} and \eqref{eq:quantization-errors} proves \eqref{eq:conditional-law-lower}.
\end{proof}

\begin{corollary}\label{cor:conditional-law-packing}
Let $\mathbb X$, $(t_0,\varepsilon_0)$, $(E,d_E)$, $R$ and $\mathcal K=\Lc(R(X)\mid\cF_{t_0}^X)$ be as in Lemma~\ref{lem:conditional-law-lower} with $p=1$.
Let $J\geq2$ and $r,w>0$, and suppose that there are laws $\mu_1,\ldots,\mu_J\in\cP_1(E)$ and pairwise disjoint events $\Omega_1,\ldots,\Omega_J\in\cF_{t_0}^X$ such that, for all $i,j\leq J$,
\begin{enumerate}[label=\textnormal{(\roman*)}]
 \item\label{it:packing-separation} $\cW_1(\mu_i,\mu_j)\geq4r$ if $i\ne j$;
 \item\label{it:packing-mass} $\P(\Omega_j)\geq w/J$;
 \item\label{it:packing-realization} $\cW_1(\mathcal K,\mu_j)\leq r$ almost surely on $\Omega_j$.
\end{enumerate}
Then, for every integer $1\leq M\leq J/2$,
\begin{equation}\label{eq:conditional-law-packing}
 e^{\mathrm{ad}}_{M,1}(\mathbb X)\geq\min\{\varepsilon_0,wr/2\}.
\end{equation}
\end{corollary}

\begin{proof}
At most $M$ of the separated laws $\mu_j$ can lie near one of $M$ centers; each of the remaining laws is realized with probability at least $w/J$ and then contributes an error of at least $r$.
To make this precise, fix an integer $1\leq M\leq J/2$ and centers $\lambda_1,\ldots,\lambda_M\in\cP_1(E)$.
By \ref{it:packing-separation}, the open $\cW_1$-balls of radius $2r$ around $\mu_1,\ldots,\mu_J$ are pairwise disjoint, so each $\lambda_i$ lies in at most one of them.
Hence at most $M$ indices $j$ satisfy $\cW_1(\mu_j,\lambda_i)<2r$ for some $i\leq M$, and the index set
\[
 \mathcal U:=\{j\leq J:\ \cW_1(\mu_j,\lambda_i)\geq2r\text{ for all }i\leq M\}
\]
satisfies $|\mathcal U|\geq J-M\geq J/2$.
For $j\in\mathcal U$ and $i\leq M$, the triangle inequality and \ref{it:packing-realization} give, almost surely on $\Omega_j$,
\[
 \cW_1(\mathcal K,\lambda_i)
 \geq \cW_1(\mu_j,\lambda_i)-\cW_1(\mathcal K,\mu_j)
 \geq2r-r=r.
\]
Since the events $\Omega_j$ are disjoint, \ref{it:packing-mass} yields
\[
 \E\min_{i\leq M}\cW_1(\mathcal K,\lambda_i)
 \geq r\sum_{j\in\mathcal U}\P(\Omega_j)
 \geq r\,|\mathcal U|\,\frac wJ
 \geq\frac{wr}2.
\]
As the centers were arbitrary, $Q_{M,1}\geq wr/2$, and Lemma~\ref{lem:conditional-law-lower} gives \eqref{eq:conditional-law-packing}.
\end{proof}

\begin{lemma}\label{lem:half-subset-packing}
There exist universal constants $c>0$ and $N_0\geq2$ such that, for every even integer $N\geq N_0$, there is a family $\mathcal A_N$ of subsets of $\{1,\ldots,N\}$ of size $N/2$ satisfying
\begin{equation}\label{eq:half-subset-code}
 |\mathcal A_N|\geq e^{cN},
 \qquad
 |A\mathbin\triangle B|\geq\frac N4\quad(A,B\in\mathcal A_N,\ A\ne B).
\end{equation}
Consequently, for every metric space $(E,d_E)$, any $r>0$ and any points $x_1,\ldots,x_N\in E$ with
\begin{equation}
 \min_{1\leq i<j\leq N}d_E(x_i,x_j)\geq r,
\end{equation}
the uniform laws $\mu_A:=(2/N)\sum_{i\in A}\delta_{x_i}$ on $\{x_i:i\in A\}$, $A\in\mathcal A_N$, satisfy
\begin{equation}
 \cW_1(\mu_A,\mu_B)\geq\frac r4,
 \qquad A,B\in\mathcal A_N,\quad A\ne B,
\end{equation}
where $\cW_1$ is computed with $d_E$.
\end{lemma}

\begin{proof}
We first construct $\mathcal A_N$, which depends only on $N$, and then derive the transport bound.
Call two subsets of $\{1,\ldots,N\}$ of size $N/2$ \emph{close} if their symmetric difference has fewer than $N/4$ elements; in particular, every such subset is close to itself.
Let $\mathcal A_N$ be a maximal family of subsets of size $N/2$, no two of which are close.
Then distinct $A,B\in\mathcal A_N$ satisfy $|A\mathbin\triangle B|\geq N/4$, and, by maximality, every subset of size $N/2$ is close to some member of $\mathcal A_N$.
To count the subsets close to a fixed $A$, note that $B\mapsto A\mathbin\triangle B$ is injective and that $3^{N/4-h}\geq1$ for $h<N/4$.
Hence at most
\[
 \sum_{h<N/4}\binom Nh
 \leq\sum_{h=0}^N\binom Nh3^{N/4-h}
 =3^{N/4}(4/3)^N
\]
subsets are close to $A$.
On the other hand, $\binom N{N/2}$ is the largest of the $N+1$ binomial coefficients $\binom Nh$, whose sum is $2^N$, so there are $\binom N{N/2}\geq2^N/(N+1)$ subsets of size $N/2$.
Every one of them is close to a member of $\mathcal A_N$, and each member is close to at most $3^{N/4}(4/3)^N$ of them; therefore
\[
 |\mathcal A_N|
 \geq\frac{\binom N{N/2}}{3^{N/4}(4/3)^N}
 \geq\frac{(3^{3/4}/2)^N}{N+1}.
\]
As $\log(3^{3/4}/2)>1/10$, the right-hand side is at least $e^{cN}$ for $c:=1/10$ and every even $N\geq N_0:=168$.
This proves \eqref{eq:half-subset-code}.

Now fix distinct $A,B\in\mathcal A_N$ and a coupling $\pi$ of $\mu_A$ and $\mu_B$.
Both laws are supported by $\{x_1,\ldots,x_N\}$, whose distinct points are at distance at least $r$.
For $i\in A\setminus B$, the point $x_i$ is not an atom of $\mu_B$, so $\pi$ moves the entire mass $2/N$ at $x_i$ over a distance of at least $r$.
Since $|A|=|B|$ implies $|A\setminus B|=|A\mathbin\triangle B|/2$, it follows that
\[
 \int d_E\,\mathrm d\pi
 \geq\frac{2r}N|A\setminus B|
 =\frac rN|A\mathbin\triangle B|
 \geq\frac r4.
\]
Taking the infimum over all couplings $\pi$ gives $\cW_1(\mu_A,\mu_B)\geq r/4$.
\end{proof}

We now turn to the amplification mechanism.
In Lemma~\ref{lem:osgood-amplifier}, we construct a scalar Osgood flow that maps the exponentially small seed $g_m(z)$ to $z$ in one unit of time.
This allows us to use a globally Lipschitz drift to generate a small seed whose amplified value is highly sensitive to the noise.
In Corollary~\ref{cor:osgood-vector-flow}, we amplify several coordinates simultaneously, verify the joint Osgood condition, and compute the comparison flow.

\begin{lemma}\label{lem:osgood-amplifier}
For every real $m>0$, define $g_m(0):=0$ and $g_m(z):=\exp(-z^{-m})$ for $z>0$.
There are $z_*>0$ and a bounded, increasing, concave function $f_m:[0,\infty)\to[0,\infty)$ satisfying
\begin{equation}
 f_m(0)=0,\qquad f_m(z)>0\quad(z>0),\qquad
 \int_{0+}\frac{\mathrm dz}{f_m(z)}=\infty,
\end{equation}
such that the flow $\Psi_t$ of $\dot z=f_m(z)$ satisfies
\begin{equation}\label{eq:amplifier-flow}
 \Psi_1(g_m(z))=z\quad(0\leq z\leq z_*),\qquad
 \Psi_1(s)=[\log(1/s)]^{-1/m}\quad(0<s\leq g_m(z_*)).
\end{equation}
For this $f_m$, the odd extension $\widetilde f_m$ also satisfies
\begin{equation}\label{eq:amplifier-modulus}
 |\widetilde f_m(a)-\widetilde f_m(b)|
 \leq2f_m(|a-b|/2),\qquad a,b\in\R.
\end{equation}
\end{lemma}

\begin{proof}
Set
\begin{equation}
 h(z):=\frac{g_m''(z)}{g_m'(z)}
 =mz^{-m-1}-\frac{m+1}{z},\qquad
 h'(z)+\frac{h(z)}z=-m^2z^{-m-2}.
\end{equation}
Choose $z_*$ sufficiently small that, on $(0,z_*]$, we have $g_m(z)<z$, $h(z)>0$, and $0<g_m'(z)\leq\vartheta<1$.
We also require $mz_*^{-m}>1$.
Put $\zeta_0:=z_*$ and $\zeta_{j+1}:=g_m(\zeta_j)$.
On $[\zeta_1,\zeta_0]$, let $f_0$ be the affine function with
\begin{equation}
 f_0(\zeta_0)=1,\qquad f_0(\zeta_1)=g_m'(\zeta_0).
\end{equation}
Its slope is positive because $g_m'(\zeta_0)<1$.
Its intercept is positive because
$\zeta_0g_m'(\zeta_0)-\zeta_1=\zeta_1(m\zeta_0^{-m}-1)>0$.
In particular, $0\leq zf_0'(z)\leq f_0(z)$ on this first interval.
Extend $f_0$ successively to $[\zeta_{j+2},\zeta_{j+1}]$ through
\begin{equation}\label{eq:amplifier-recursion}
 f_0(g_m(z))=g_m'(z)f_0(z),\qquad 0<z\leq z_*.
\end{equation}
The chosen endpoint values make this extension continuous at every $\zeta_j$.

We verify monotonicity and concavity together with $zf_0'(z)\leq f_0(z)$ on each smooth interval.
Differentiating \eqref{eq:amplifier-recursion} gives
\begin{equation}
 f_0'(g_m(z))=f_0'(z)+h(z)f_0(z),
 \qquad
 f_0''(g_m(z))=
 \frac{f_0''(z)+h(z)f_0'(z)+h'(z)f_0(z)}{g_m'(z)}.
\end{equation}
If these properties hold at $z$, the first identity gives a nonnegative derivative at $g_m(z)$.
The numerator in the second identity is at most
$(h(z)/z+h'(z))f_0(z)<0$.
Moreover,
\begin{equation}
 f_0(g_m(z))-g_m(z)f_0'(g_m(z))
 =g_m(z)\left[\frac{m+1}{z}f_0(z)-f_0'(z)\right]\geq0.
\end{equation}
The last inequality uses $zf_0'(z)\leq f_0(z)$ and $m>0$.
This proves the three properties by induction within every interval.
At $\zeta_1$, the left derivative exceeds the right derivative by $h(\zeta_0)>0$.
Thus the derivative decreases across the first junction.
Taking one-sided limits in the derivative recursion and using the continuity of $h$ and $f_0$, we obtain, for $j\geq1$,
\begin{equation}
 f_0'(\zeta_{j+1}-)-f_0'(\zeta_{j+1}+)
 =f_0'(\zeta_j-)-f_0'(\zeta_j+)
 =h(\zeta_0)>0.
\end{equation}
Consequently, $f_0$ is increasing and concave on $(0,z_*]$.
On $[\zeta_{j+1},\zeta_j]$, $j\geq0$, its supremum is at most $\vartheta^j$.
Setting $f_0(0):=0$ therefore extends it continuously and concavely to $[0,z_*]$.

For $0<z\leq z_*$, define $I(z):=\int_{g_m(z)}^z\mathrm du/f_0(u)$.
Continuity and positivity of $f_0$ give
\begin{equation}
 I'(z)=\frac1{f_0(z)}-\frac{g_m'(z)}{f_0(g_m(z))}=0.
\end{equation}
Thus $I(z)$ is a fixed number $I\in(0,\infty)$.
Set $f_m:=If_0$ on $[0,z_*]$ and extend it constantly on $[z_*,\infty)$.
This extension is bounded, increasing, and concave.
Every interval $[\zeta_{j+1},\zeta_j]$ has integral $1$ against $1/f_m$, which proves the Osgood divergence.
The same relation also proves \eqref{eq:amplifier-flow} via the change of variable formula
\begin{equation}
\int_s^{\Psi_t(s)}\frac{du}{f_m(u)} = t.
\end{equation}
Finally, concavity and $f_m(0)=0$ imply subadditivity.
For $a,b$ of the same sign, $|\widetilde f_m(a)-\widetilde f_m(b)|\leq f_m(|a-b|)\leq2f_m(|a-b|/2)$.
For opposite signs, \eqref{eq:amplifier-modulus} follows from Jensen's inequality.
\end{proof}

\begin{corollary}\label{cor:osgood-vector-flow}
Fix a real $m>0$ and an integer $q\geq1$, and let $f_m,\widetilde f_m,\Psi$ be as in Lemma~\ref{lem:osgood-amplifier}.
Set $F_q(z):=(\widetilde f_m(z_1),\ldots,\widetilde f_m(z_q))$ and
\begin{equation}\label{eq:example-joint-modulus}
 \rho_q(s):=4\sqrt q\,\sqrt s\,
 f_m\!\left(\frac{\sqrt s}{2\sqrt q}\right),\qquad s\geq0.
\end{equation}
Let $S:\R^2\to\R^q$ be globally Lipschitz with constant $L_S$.
Suppose $\kappa,\chi$ are nonnegative continuous functions on $[0,1]$, each with integral one, and $\sup\supp\kappa<\inf\supp\chi$.
Then the drift
$b(t,u,v,z):=(0,0,\kappa(t)S(u,v)+\chi(t)F_q(z))$,
together with any state-independent diffusion coefficient, satisfies Assumption~\ref{ass:osgood} with $\ell=L_S\kappa$, $k=\chi$, and $\rho=\rho_q$.
Its comparison flow satisfies, for sufficiently small $\xi>0$,
\begin{equation}\label{eq:example-comparison-flow}
 \Phi(\xi)
 =4q\Psi_1\!\left(\frac{e^{L_S/2}\sqrt\xi}{2\sqrt q}\right)^2
 \asymp[\log(1/\xi)]^{-2/m}.
\end{equation}
In particular, Assumption~\ref{ass:flow-rate} holds with $\beta=2/m$.
\end{corollary}

\begin{proof}
The function $\rho_q$ is increasing, vanishes at zero, and is positive elsewhere.
On each smooth piece, its concavity is equivalent to
$r f_m'(r)+r^2f_m''(r)-f_m(r)\leq0$, which follows from concavity of $f_m$ and $f_m(0)=0$.
Its derivative decreases at the junctions, and continuity extends concavity to zero.
The substitution $s=4qr^2$ gives
$\int_{0+}\mathrm ds/\rho_q(s)=\int_{0+}\mathrm dr/f_m(r)=\infty$.
The same concavity calculation, \eqref{eq:amplifier-modulus}, and Jensen's inequality give
\begin{align}
 2\langle z-z',F_q(z)-F_q(z')\rangle
 &\leq4\sum_{i=1}^q |z_i-z_i'|f_m(|z_i-z_i'|/2)\nonumber\\
 &\leq4\sqrt q\,|z-z'|f_m\!\left(\frac{|z-z'|}{2\sqrt q}\right)
 =\rho_q(|z-z'|^2).
\end{align}
For $x=(u,v,z)$ and $x'=(u',v',z')$, the Lipschitz continuity of $S$ gives
$2\langle z-z',S(u,v)-S(u',v')\rangle\leq L_S|x-x'|^2$.
Since $\rho_q$ is increasing and the diffusion coefficient is state-independent, these estimates prove \eqref{eq:joint-osgood}.
Finally, denoting by $z(t;\xi)=4qr(t)^2$ the correspongding solution of \eqref{eq:comparison-flow} gives
$\dot r=(L_S/2)\kappa(t)r+\chi(t)f_m(r)$ with $r(0)=\sqrt\xi/(2\sqrt q)$.
Since $\sup\supp\kappa<\inf\supp\chi$ and both functions have integral one, we have $r(1)=\Psi_1(e^{L_S/2}r(0))$, which is the first identity in \eqref{eq:example-comparison-flow}.
The second equality in \eqref{eq:example-comparison-flow} follows directly from \eqref{eq:amplifier-flow} 
\end{proof}

\subsection{Proof of results in Section \ref{sec:notation}}
\label{sec:notation-proofs}

\begin{proof}[Proof of Proposition~\ref{prop:iterated-log-flow}]
Choose $u_*>0$ so small that \eqref{eq:iterated-log-modulus} holds and all its logarithmic factors are at least one for $0<u\leq u_*$.
Then
\begin{equation}
 \int_a^{u_*}\frac{\mathrm du}{\rho(u)}
 \geq\frac1C\left[\log^{(j+1)}(1/a)-\log^{(j+1)}(1/u_*)\right]
 \longrightarrow\infty\qquad\text{as }a\downarrow0,
\end{equation}
which proves the Osgood condition.
Boundedness of $\ell,k$ also gives a constant $C_0>0$ such that
\begin{equation}
 \ell(t)u+k(t)\rho(u)\leq C_0u\prod_{i=1}^j\log^{(i)}(1/u),
 \qquad t\in[0,1],\quad 0<u\leq u_*.
\end{equation}

Fix $\beta>0$ and write $a=e^{-s}$.
For large $s$, we have $e^{-s}<s^{-\beta}<u_*$.
The solution of $w'=C_0w\prod_{i=1}^j\log^{(i)}(1/w)$ with $w(0)=e^{-s}$ reaches $s^{-\beta}$ at time
\begin{equation}
 \tau_s=\frac1{C_0}\int_{e^{-s}}^{s^{-\beta}}
 \frac{\mathrm du}{u\prod_{i=1}^j\log^{(i)}(1/u)}
 =\frac1{C_0}\left[\log^{(j)}s-\log^{(j)}(\beta\log s)\right].
\end{equation}
For fixed $j$ and $\beta$, we have $\log^{(j)}(\beta\log s)=\log^{(j+1)}s+O(1)$, so $\tau_s\to\infty$.
Scalar comparison applies on $[0,1\wedge\tau_s]$, since the comparison trajectory remains in $(0,u_*]$.
Thus, for large $s$ with $\tau_s>1$,
\begin{equation}
 \Phi(a)=z(1;a)\leq w(1)<s^{-\beta}=[\log(1/a)]^{-\beta}.
\end{equation}
This proves \eqref{eq:flow-condition} for every $\beta>0$.
The choice $\rho(u)=u$ satisfies \eqref{eq:iterated-log-modulus} with $j=1$.
\end{proof}

\begin{proof}[Proof of Proposition~\ref{prop:strict-finite-support}]
Let $S:=\supp\nu$.
Every point of $S$ has positive $\nu$-mass, so completion does not change the restriction of the canonical filtration to $S$.
On this finite set, the $\sigma$-fields $\cF_t^\nu|_S$ form an increasing chain, and each strict increase splits at least one atom.
There are therefore at most $|S|-1$ strict increases.
Right-continuity places each change at the left endpoint of the next interval, so there is a deterministic partition $0=t_0<\cdots<t_N=1$ such that, modulo $\nu$-null sets,
\begin{equation}
 \cF_t^\nu=\cF_{t_k}^\nu,\qquad t\in[t_k,t_{k+1}),\quad 0\leq k<N.
\end{equation}
The usual Brownian filtration is trivial at time zero by Blumenthal's zero--one law and satisfies $\cF_{t-}^{\W_d}=\cF_t^{\W_d}$ for every $t\in(0,1]$.
These are precisely the filtration conditions in \citet[Proposition~C.2]{BBPSZ2025}, which implies that the coordinate processes are independent under every strictly bicausal coupling.
The product coupling is itself strictly bicausal, which proves the equality in \eqref{eq:strict-finite-support-cost}.

The function $g(y):=\E\|B-y\|_\infty^p$ is convex, and the symmetry of Brownian motion gives $g(-y)=g(y)$.
Therefore $g(0)\leq\tfrac12(g(y)+g(-y))=g(y)$, proving the lower bound in \eqref{eq:strict-finite-support-cost}.
Taking $\nu=\delta_0$, where $0$ denotes the zero path, attains this bound for every $M\geq1$ and proves \eqref{eq:strict-quantization-obstruction}.
\end{proof}

\subsection{Proof of Proposition~\ref{prop:sampling-sharp}}
\label{sec:raw-empirical-proof}

\begin{proof}
For every $i\ne j$ and $n\in\N$, the vector $B_{1/n}^i-B_{1/n}^j$ has a nondegenerate Gaussian law.
Hence, on a single probability-one event, $B_{1/n}^i\ne B_{1/n}^j$ for all such $i,j,n$.
Fix a realization in this event and $M\in\N$, and write $\nu$ for the realized empirical law $\mu_M$.
Under $\nu$, the value of the canonical path at time $1/n$ identifies its support point and hence its entire trajectory.
Right-continuity therefore gives $\cF_0^\nu=\bigcap_{n\geq1}\cF_{1/n}^\nu=\cF_1^\nu$ modulo $\nu$-null sets.

Fix $p\in[1,\infty)$, $\varepsilon\in[0,1)$, and $\pi\in\Pi_{bc}^{\varepsilon}(\W_d,\Can(\nu))$.
Let $X$ and $Y$ denote the Brownian and empirical coordinate processes under $\pi$.
By Definition~\ref{def:bicausal-couplings} at time zero and the independent increments of the Brownian marginal,
\begin{equation}
 \Law_\pi\bigl(X_1-X_\varepsilon\mid\cF_\varepsilon^X\vee\cF_0^Y\bigr)
 =\Law_\pi\bigl(X_1-X_\varepsilon\mid\cF_\varepsilon^X\bigr)
 =\mathcal N\bigl(0,(1-\varepsilon)I_d\bigr).
\end{equation}
Since $Y_1$ is $\cF_0^Y$-measurable, the shift $X_\varepsilon-Y_1$ is measurable with respect to $\cF_\varepsilon^X\vee\cF_0^Y$.
By symmetry and convexity, the function $a\mapsto\E|a+\sqrt{1-\varepsilon}\,G|^p$ is minimized at $a=0$.
Consequently,
\begin{equation}
 \E_\pi\bigl[|X_1-Y_1|^p\mid\cF_\varepsilon^X\vee\cF_0^Y\bigr]
 \geq(1-\varepsilon)^{p/2}m_{p,d}^p.
\end{equation}
Using $\|X-Y\|_\infty\geq|X_1-Y_1|$, we obtain
\begin{equation}
 \varepsilon+\bigl(\E_\pi\|X-Y\|_\infty^p\bigr)^{1/p}
 \geq\varepsilon+m_{p,d}\sqrt{1-\varepsilon}
 \geq\min\{1,m_{p,d}\}.
\end{equation}
The last inequality follows from $\sqrt{1-\varepsilon}\geq1-\varepsilon$.
Delays $\varepsilon\geq1$ already contribute at least one to the cost.
Taking both infima in \eqref{eq:aw-definition} proves \eqref{eq:raw-empirical-lower}.
Taking expectation over the samples gives the same lower bound for the mean error.
\end{proof}

\subsection{Proof of Theorem~\ref{thm:recursive-sampling}}
\label{sec:main-proofs}

\begin{proof}[Proof of Theorem~\ref{thm:recursive-sampling}]
Fix $p\geq1$ and a growth bound $K\geq1$ as in Assumption~\ref{ass:coefficients}.
Write $\mu:=\Lc(X)$, $h:=1/n$, and $t_k:=kh$.
Initially, take $h,r>0$ small enough for Lemma~\ref{lem:local-approximation}, $0<\delta\leq1/n$, and $R\geq1$.
Constants below may depend on the fixed data and on the indicated moment order, but not on $n,L,r,\delta,R$ or the sampling realization.

\textbf{Step 1: Constructing the coupling.}
We first define the labels used to match diffusion segments.
For a segment $s\in C([t_k,t_{k+1}];\R^d)$ starting at $y$, put
\begin{equation}
 V_{k,y}(s)(u):=\frac{s(t_k+hu)-y}{\sqrt h(1+|y|)},\qquad 0\leq u\leq1.
\end{equation}
Use the codebook $\cA_r$ and encoder $\alpha_{h,r}$ from Lemma~\ref{lem:local-approximation}.
Assign the bad label $\perp$ if
$\|V_{k,y}(s)-\alpha_{h,r}(V_{k,y}(s))\|_\infty>C_dr$ or $s(t_{k+1})\notin[-R,R]^d$.
Otherwise, assign the good label consisting of the codeword $\alpha_{h,r}(V_{k,y}(s))$ and the endpoint cell
\begin{equation}
 \begin{cases}
 \left\lfloor f_{t_{k+1}}(s(t_{k+1}))/\delta\right\rfloor,&d=1,\\
 \left(\left\lfloor s_i(t_{k+1})/\delta\right\rfloor\right)_{i=1}^d,&d\geq2,
 \end{cases}
\end{equation}
where $f_t$ is defined in Lemma~\ref{lem:scalar-endpoints}.
Denote this Borel label map by $\Lambda_{k,y}$.
In one dimension, the growth bound for $f_t$ in that lemma bounds the length of $f_t([-R,R])$ by $C(1+R)$.
In higher dimensions, each coordinate has at most $C(1+R)/\delta$ possible cells.
Together with \eqref{eq:forest-size}, these observations give a common bound $A$ on the number of labels, including $\perp$, with
\begin{equation}\label{eq:sampling-label-count}
 A\leq1+C\exp(C_0/r^2)\left(\frac{1+R}{\delta}\right)^d.
\end{equation}
If two segments $s,s'$ start at the same $y$ and have the same good label, the triangle inequality and the endpoint cells give
\begin{align}
 &\|s-s'\|_{\infty,[t_k,t_{k+1}]}\leq C\sqrt h(1+|y|)r,\notag\\
 &|f_{t_{k+1}}(s(t_{k+1}))-f_{t_{k+1}}(s'(t_{k+1}))|\leq\delta\quad(d=1),\notag\\
 &|s(t_{k+1})-s'(t_{k+1})|\leq\sqrt d\,\delta\quad(d\geq2).
 \label{eq:sampling-label-errors}
\end{align}

Fix a sampling realization $\mathcal S$ and set $\nu_{\mathcal S}:=\mu^{n,L}(\mathcal S)$.
For each distinct sampled prefix $g$ on $[0,t_k]$, let
\begin{equation}
 I_k(g):=\{\mathbf i\in\mathcal I_k:g_{\mathbf i}=g\}.
\end{equation}
Each index in $\mathcal I_k$ has $L^{n-k}$ descendants, so the prefix $g$ has $\nu_{\mathcal S}$-mass $|I_k(g)|L^{-k}$.
Conditional on this entire prefix, the next segment has law
\begin{equation}\label{eq:sampling-prefix-kernel}
 Q^{\mathcal S}_{k,g}:=
 \frac1{|I_k(g)|L}\sum_{\mathbf i\in I_k(g)}\sum_{j=1}^L\delta_{S_{\mathbf i,j}}.
\end{equation}

On an auxiliary space, take a Brownian motion $B$ and independent uniform variables $U_0,\ldots,U_{n-1}$, all independent of the sampling.
Let $\mathcal H$ be the usual filtration in the proof of Lemma~\ref{lem:block-bicausal}, with $U_k$ revealed at $t_{k+1}$, and use $B$ to drive $X$.
Starting from $Y_0=x_0$, suppose that $Y|_{[0,t_k]}=g$.
At every step, we solve the SDE from $g(t_k)=Y_{t_k}$ on $[t_k,t_{k+1}]$ with the same Brownian increment as $X$, obtaining the comparison segment $C^{(k)}$.
Conditionally on $\mathcal H_{t_k}$, its law is $P_k(g(t_k),\cdot)$.
Let $\lambda_{k,g(t_k)}$ be its label law under $\Lambda_{k,g(t_k)}$, and let $\widehat\lambda_{k,g}$ be the label law of \eqref{eq:sampling-prefix-kernel}.
Lemma~\ref{lem:maximal-label-coupling}\textnormal{(i)}, using $U_k$, selects an entry of this pool with conditional law $Q^{\mathcal S}_{k,g}$ and mismatch probability $\TV(\lambda_{k,g(t_k)},\widehat\lambda_{k,g})$.
Use the selected segment as $Y|_{[t_k,t_{k+1}]}$.
At the next step, recompute the pool from the observed prefix alone.
Induction gives $\Lc(X)=\mu$ and $\Lc(Y)=\nu_{\mathcal S}$.

By the Borel strong-solution maps and Lemma~\ref{lem:maximal-label-coupling}, the coupling construction and its cost are Borel measurable in $\mathcal S$.
For fixed $\mathcal S$, denote the resulting segment map, rescaled to $[0,1]$, by $F_k^{\mathcal S}(g,w,u)$.
For prefixes $g$ outside the sampled set, set $F_k^{\mathcal S}(g,w,u)(v):=g(t_k)$ for $0\leq v\leq1$.
These maps are Borel and satisfy $F_k^{\mathcal S}(g,w,u)(0)=g(t_k)$, so Lemma~\ref{lem:block-bicausal} gives $h$-bicausality.
Writing $\E^{\mathcal S}$ for expectation over the auxiliary randomness, \eqref{eq:aw-definition} yields
\begin{equation}\label{eq:sampling-bicausal-cost}
 \AW_p(\mathbb X,\Can(\nu_{\mathcal S}))
 \leq\frac1n+\left(\E^{\mathcal S}\|X-Y\|_\infty^p\right)^{1/p}.
\end{equation}

\textbf{Step 2: Error estimate under $\mathcal S$.}
Let $E_k$ be the event that $C^{(k)}$ and $Y|_{[t_k,t_{k+1}]}$ have the same good label, and put
\begin{equation}
 G_k:=\bigcap_{j<k}E_j,\qquad G_0:=\Omega.
\end{equation}
Construct an $\mathcal H$-adapted process $Z^*$ from $x_0$ by solving the original SDE with $B$ between grid times.
As long as all preceding blocks have succeeded, its initial state on block $k$ is $Y_{t_k}$, so it follows $C^{(k)}$.
At the end of a successful block, reset it to $Y_{t_{k+1}}$.
At the first failure, make no reset and thereafter continue from its own state without further resets.
Equivalently, its jumps are
\begin{equation}
 \Delta Z^*_{t_{k+1}}
 =\1_{G_{k+1}}\bigl(Y_{t_{k+1}}-C^{(k)}_{t_{k+1}}\bigr),\qquad 0\leq k<n.
\end{equation}
The event $G_{k+1}$ and both endpoints are $\mathcal H_{t_{k+1}}$-measurable.
Thus these jumps are adapted, including the possible jump at $t_n=1$.
On $G_{k+1}$, \eqref{eq:sampling-label-errors} gives the transformed jump bound \eqref{eq:scalar-adapted-jumps} in one dimension and the bound $|\Delta Z^*_{t_{k+1}}|\leq\sqrt d\,\delta$ in higher dimensions.
All other jumps vanish.
Applying Lemma~\ref{lem:scalar-endpoints} when $d=1$ and Lemma~\ref{lem:osgood-endpoints} under Assumption~\ref{ass:osgood} when $d\geq2$, we obtain
\begin{equation}\label{eq:sampling-fixed-stability}
 \left(\E^{\mathcal S}\|X-Z^*\|_\infty^p\right)^{1/p}
 \leq C_p\mathsf S_p(n\delta),
\end{equation}
uniformly in $\mathcal S$, where
\begin{equation}\label{eq:strategy-stability}
 \mathsf S_p(a):=
 \begin{cases}
 a^{1/[5(2p-1)]},&d=1,\\
 \Phi(C_pa)^{1/(4p-2)},&d\geq2.
 \end{cases}
\end{equation}
On $G_n$, we have $Z^*_{t_k}=Y_{t_k}$ at every grid time.
Applying the path estimate in \eqref{eq:sampling-label-errors} on each block gives
\begin{equation}\label{eq:sampling-good-error}
 \|Z^*-Y\|_\infty\1_{G_n}
 \leq C\frac r{\sqrt n}(1+\|Y\|_\infty)\1_{G_n}.
\end{equation}

\textbf{Step 3: Averaging over sampling.}
We now write $\E:=\E_{\mathcal S}\E^{\mathcal S}$ and $\P$ for the joint expectation and probability.
For each deterministic full index $\mathbf i\in\mathcal I_n$, its successive branches are sampled from the transition laws $P_k$.
Induction over $k$, using the Markov property of the diffusion, shows that $g_{\mathbf i}$ has law $\mu$ after averaging over sampling.
Consequently, for every nonnegative Borel function $F$ on path space,
\begin{equation}\label{eq:sampling-average-law}
 \E_{\mathcal S}\int F(y)\,\mu^{n,L}(\mathrm dy)
 =\frac1{L^n}\sum_{\mathbf i\in\mathcal I_n}\E_{\mathcal S}F(g_{\mathbf i})
 =\int F(x)\,\mu(\mathrm dx).
\end{equation}
In particular, $X$ and $Y$ both have law $\mu$ under the joint probability, so their maximal moments of every finite order are bounded uniformly in $n,L$.

Let $\mathcal D_k$ be the $\sigma$-field generated by all sampled branches up to time $t_k$.
For fixed $\mathcal S$, the mismatch probability on block $k$ is
\begin{equation}
 \sum_g\frac{|I_k(g)|}{L^k}
 \TV(\lambda_{k,g(t_k)},\widehat\lambda_{k,g}),
\end{equation}
where the sum runs over distinct sampled prefixes.
Conditionally on $\mathcal D_k$, these prefixes and their multiplicities are fixed, and the $|I_k(g)|L$ entries in each pool are i.i.d.\ with law $P_k(g(t_k),\cdot)$.
The empirical estimate in Lemma~\ref{lem:maximal-label-coupling}\textnormal{(ii)}, applied to each such pool, gives
\begin{equation}
 \E_{\mathcal S}\left[
 \sum_g\frac{|I_k(g)|}{L^k}\TV(\lambda_{k,g(t_k)},\widehat\lambda_{k,g})
 \,\middle|\,\mathcal D_k\right]
 \leq\frac12\sum_g\frac{|I_k(g)|}{L^k}
 \sqrt{\frac A{|I_k(g)|L}}
 \leq\frac12\sqrt{\frac A L}.
\end{equation}
Here $\sum_g|I_k(g)|L^{-k}=1$ and $|I_k(g)|\geq1$.
Taking expectations bounds the joint mismatch probability by $\frac12\sqrt{A/L}$ on every block.

For each fixed $\mathcal S$, apply Lemma~\ref{lem:local-approximation} to $C^{(k)}_{t_k+u}$, $0\leq u\leq h$, with the shifted filtration $(\mathcal H_{t_k+u})_{u\leq h}$ and Brownian increment $B_{t_k+u}-B_{t_k}$.
Its initial state $Y_{t_k}$ is $\mathcal H_{t_k}$-measurable, so the coding failure probability is at most $e^{-c/r^2}+Ce^{-cn}$ conditionally on $\mathcal H_{t_k}$.
Averaging this bound over $\mathcal S$ gives the same bound under the joint probability.
Moreover, \eqref{eq:sampling-average-law} gives $\Lc(Y_{t_k})=\Lc(X_{t_k})$ under the joint probability.
Since $C^{(k)}$ has conditional law $P_k(Y_{t_k},\cdot)$, its endpoint has the same law as $X_{t_{k+1}}$ after averaging.
For every finite $q\geq1$, Markov's inequality therefore yields
\begin{equation}
 \P\bigl(\Lambda_{k,Y_{t_k}}(C^{(k)})=\perp\bigr)
 \leq e^{-c/r^2}+Ce^{-cn}+C_qR^{-q}.
\end{equation}
Failure of $E_k$ requires either a bad comparison label or a mismatch.
A union bound over the $n$ blocks along the selected path gives
\begin{equation}\label{eq:sampling-failure-probability}
 \P(G_n^c)\leq\min\{1,Q\},\qquad
 Q:=n\left[e^{-c/r^2}+Ce^{-cn}+C_qR^{-q}+\frac12\sqrt{\frac A L}\right].
\end{equation}

On $G_n$, use \eqref{eq:sampling-fixed-stability}, \eqref{eq:sampling-good-error}, and the moment bounds from \eqref{eq:sampling-average-law}.
On $G_n^c$, H\"older's inequality gives
\begin{equation}
 \E\bigl[\|X-Y\|_\infty^p\1_{G_n^c}\bigr]
 \leq\bigl(\E\|X-Y\|_\infty^{2p}\bigr)^{1/2}\P(G_n^c)^{1/2}
 \leq C_p\P(G_n^c)^{1/2}.
\end{equation}
Combining the two events yields
\begin{equation}
 \bigl(\E\|X-Y\|_\infty^p\bigr)^{1/p}
 \leq C_p\left[\mathsf S_p(n\delta)+\frac r{\sqrt n}+\P(G_n^c)^{1/(2p)}\right].
\end{equation}
Finally, average \eqref{eq:sampling-bicausal-cost} over $\mathcal S$ and use Jensen's inequality and \eqref{eq:sampling-failure-probability} to obtain
\begin{equation}\label{eq:sampling-master-bound}
 \mathcal E_{n,L,p}\leq\frac1n+
 C_p\left[\mathsf S_p(n\delta)+\frac r{\sqrt n}+\min\{1,Q\}^{1/(2p)}\right].
\end{equation}

\textbf{Step 4: Choosing the scales.}
In one dimension, take $\delta:=n^{-(10p+1)}$, $R:=n$, $q:=2p+1$, and $r:=(a\log L)^{-1/2}$, where $a>0$ is fixed and small enough that $C_0a\leq1/4$ in \eqref{eq:sampling-label-count}.
Then
\begin{equation}
 \mathsf S_p(n\delta)=n^{-2p/(2p-1)}\leq n^{-1},\qquad
 A\leq CL^{1/4}n^{10p+2}.
\end{equation}
Choose $a_p\geq\max\{4(10p+2),\,4(2p+1),\,(2p+1)/(ca)\}$, where $c$ is the constant in $e^{-c/r^2}$.
For $L\geq n^{a_p}$, substituting these choices and the bound on $A$ into \eqref{eq:sampling-failure-probability} gives $Q\leq C_pn^{-2p}$.
Substitution into \eqref{eq:sampling-master-bound}, with $r/\sqrt n=a^{-1/2}(n\log L)^{-1/2}$, proves \eqref{eq:scalar-sampling-rate}.

In higher dimensions, set $\delta:=L^{-1/(8d)}/(C_*n)$, $R:=n$, $q:=2p+1$, and $r:=(a\log L)^{-1/2}$.
Choose $C_*\geq1$ to dominate the constant inside $\Phi$ in the definition of $\mathsf S_p$ in \eqref{eq:strategy-stability}, and choose $a>0$ so that $C_0a\leq1/8$.
Monotonicity of $\Phi$ and \eqref{eq:sampling-label-count} give
\begin{equation}
 \mathsf S_p(n\delta)\leq\Phi(L^{-1/(8d)})^{1/(4p-2)},\qquad
 A\leq CL^{1/4}n^{2d}.
\end{equation}
Choose $a_p\geq\max\{8d,\,4(2p+1),\,(2p+1)/(ca)\}$.
For $L\geq n^{a_p}$, substitution into \eqref{eq:sampling-failure-probability} again gives $Q\leq C_pn^{-2p}$.
Consequently, \eqref{eq:sampling-master-bound} yields
\begin{equation}
 \mathcal E_{n,L,p}\leq C_p\left[
 \frac1n+\frac1{\sqrt{n\log L}}+\Phi(L^{-1/(8d)})^{1/(4p-2)}\right].
\end{equation}
By Assumption~\ref{ass:flow-rate}, for all sufficiently large $L$,
\begin{equation}
 \Phi(L^{-1/(8d)})^{1/(4p-2)}\leq C_p(\log L)^{-\beta/(4p-2)},
\end{equation}
which proves \eqref{eq:multidimensional-sampling-rate}.
In each case, increasing $n_p$ ensures that $h=1/n$ and $r=(a\log L)^{-1/2}$ satisfy the thresholds of Lemma~\ref{lem:local-approximation} for every $n\geq n_p$ and $L\geq\max\{2,n^{a_p}\}$.
In the higher-dimensional case it also ensures the use of Assumption~\ref{ass:flow-rate}.
\end{proof}

\subsection{Proof of Proposition~\ref{prop:counterexample}}
\label{sec:counterexample-proof}

\begin{proof}[Proof of Proposition~\ref{prop:counterexample}]
Fix $\beta>0$ and set $m:=2/\beta$ and $q:=\lceil m^2\rceil$.
Put $d:=q+2$ and let $f_m,g_m,z_*$ be as in Lemma~\ref{lem:osgood-amplifier}.
The construction encodes many terminal laws in globally Lipschitz seeds and recovers their scale through $\Psi_1(g_m(z))=z$.
Following Section~\ref{sec:proof-strategy}, we construct the separated laws, the seeds, and the diffusion in Steps~1--3, choosing the scales as their constraints arise, and derive the upper and lower bounds in Step~4.

\textbf{Step 1: Separated laws and continuous samplers.}
Let $a_k\downarrow0$ be a sequence to be chosen in Step~2.
At level $k\geq k_0\geq1$, place $2^k$ equally spaced points on each axis, with spacing $r_k>0$, to form the grid
\begin{equation}
 G_k:=\{(a_k+v_1r_k,\ldots,a_k+v_qr_k):
 v_i\in\{0,\ldots,2^k-1\}\}.
\end{equation}
We take $r_k:=a_k2^{-k}$, so that $G_k\subset[a_k,2a_k]^q$ consists of $2^{qk}$ points with minimum distance $r_k$.
Lemma~\ref{lem:half-subset-packing} provides a constant $c_*>0$ such that, for sufficiently large $k_0$, there are at least $e^{c_*2^{qk}}$ separated uniform laws on half-subsets of $G_k$.
We choose $J_k:=\lfloor e^{c_*2^{qk}}\rfloor$ of these laws and denote them by $\mu_{k,j}$.
They satisfy
\begin{equation}\label{eq:counterexample-law-separation}
 \cW_1(\mu_{k,j},\mu_{k,j'})\geq r_k/4\qquad(j\ne j').
\end{equation}
Represent each law by a step sampler on $2^{qk-1}$ equal subintervals of $[0,1]$.
Replace its jumps by affine transitions of width $\delta_k:=c_q2^{-(q+1)k}$.
For sufficiently small $c_q>0$, these intervals are disjoint, and the resulting continuous map $H_{k,j}:[0,1]\to[a_k,2a_k]^q$ satisfies, for $R_0$ uniform on $(0,1)$,
\begin{equation}\label{eq:counterexample-continuous-laws}
 \cW_1\bigl(\Lc(H_{k,j}(R_0)),\mu_{k,j}\bigr)\leq r_k/100,
 \qquad
 \operatorname{Lip}(H_{k,j})\leq C_qa_k2^{(q+1)k}.
\end{equation}
Indeed, the changes occupy length at most $2^{qk}\delta_k=c_q2^{-k}$ and have size at most $\sqrt q\,a_k$, while each transition has slope at most $\sqrt q\,a_k/\delta_k$.

\textbf{Step 2: Lipschitz seeds.}
To select among these laws, assign level $k$ the interval $I_k:=(2^{-k},2^{1-k})$ of length $2^{-k}$.
These intervals are disjoint and contained in $[0,1]$.
Divide $I_k$ into $J_k$ equal slots.
On the $j$th slot, let $\theta_{k,j}$ be continuous and piecewise affine, equal to one on the middle half and zero outside the slot, with $0\leq\theta_{k,j}\leq1$ and $\operatorname{Lip}(\theta_{k,j})\leq4J_k/|I_k|$.
Writing $\Phi_{\mathrm G}$ for the standard Gaussian distribution function, define
\begin{equation}\label{eq:example-seeds}
 T(u,v):=\sum_{k\geq k_0}\sum_{j=1}^{J_k}
 \theta_{k,j}(u)H_{k,j}(\Phi_{\mathrm G}(v)),\qquad
 S_i(u,v):=g_m(T_i(u,v)),\quad1\leq i\leq q.
\end{equation}
At most one summand is nonzero.
The selectors vanish at every slot boundary, and $a_k\to0$ gives continuity at the accumulation point $u=0$.
On each layer,
\begin{equation}
 \operatorname{Lip}(T|_{I_k\times\R})
 \leq C_qa_k\left(2^kJ_k+2^{(q+1)k}\right).
\end{equation}
Increase $k_0$ so that $2a_k\leq z_*$ and $g_m'$ is increasing on $[0,2a_k]$.
Then
\[
 \sup_{0\leq z\leq2a_k}g_m'(z)
 \leq C_m a_k^{-m-1}e^{-(2a_k)^{-m}}.
\]
These estimates give a uniform Lipschitz bound for $S$ on all layers whenever
\begin{equation}\label{eq:example-seed-condition}
 \sup_{k\geq k_0}a_k^{-m}e^{-(2a_k)^{-m}}
 \left(2^kJ_k+2^{(q+1)k}\right)<\infty.
\end{equation}
To satisfy \eqref{eq:example-seed-condition}, we choose
\begin{equation}
 a_k:=A2^{-qk/m},\qquad (2A)^{-m}>c_*,
\end{equation}
with $A>0$ fixed, and increase $k_0$ as needed.
Then $e^{-(2a_k)^{-m}}=e^{-(2A)^{-m}2^{qk}}$, whereas $J_k\leq e^{c_*2^{qk}}$ and all remaining factors grow at most exponentially in $k$.
The strict inequality $(2A)^{-m}>c_*$ therefore verifies \eqref{eq:example-seed-condition}.
For fixed $v$, the bound extends across slots and layers because $S$ vanishes at their boundaries and is continuous at zero.
For fixed $u$, the bound in $v$ is uniform.
Thus $S:\R^2\to\R^q$ is globally Lipschitz.

\textbf{Step 3: The diffusion and its conditional laws.}
Choose smooth functions $h_U,h_V$ supported in $(0,1/8)$ and $(1/2,5/8)$, respectively, with
$\int_0^1h_U^2=\int_0^1h_V^2=1$.
Choose nonnegative smooth functions $\kappa,\chi$ supported in $(11/16,3/4)$ and $(13/16,15/16)$, respectively, each with integral one.
For $F_q(z):=(\widetilde f_m(z_1),\ldots,\widetilde f_m(z_q))$, define $X=(U,V,Z)$ from zero by
\begin{equation}\label{eq:counterexample-sde}
 \mathrm dU_t=h_U(t)\,\mathrm dB_t,\qquad
 \mathrm dV_t=h_V(t)\,\mathrm dB_t,\qquad
 \mathrm dZ_t=[\kappa(t)S(U_t,V_t)+\chi(t)F_q(Z_t)]\,\mathrm dt,
\end{equation}
using a scalar Brownian motion $B$.
The drift is continuous and bounded, and the diffusion coefficient is smooth, bounded, and independent of the state.
From every deterministic starting pair $(s,x)$, the first two coordinates have explicit stochastic integral solutions.
Since $\kappa$ and $\chi$ have disjoint, ordered supports, we solve (on different time intervals) $Z$ by integrating $\kappa(t)S(U_t,V_t)$ and then by solving $\dot Z_t^i=\chi(t)\widetilde f_m(Z_t^i)$ coordinatewise, using a deterministic time change of the scalar flow.
Lemma~\ref{lem:osgood-amplifier} and \eqref{eq:amplifier-modulus} give unique global solutions for the latter, verifying Assumption~\ref{ass:coefficients}.

For the solution from zero, $U_t=U_1$ and $V_t=V_1$ for $t\geq5/8$.
Since $\int_0^1\kappa(t)\,\mathrm dt=1$, integration gives $Z_{3/4}=S(U_1,V_1)$.
The subsequent evolution, together with $\int_0^1\chi(t)\,\mathrm dt=1$ and \eqref{eq:amplifier-flow}, yields
\begin{equation}\label{eq:example-terminal-output}
 Z_1=\Psi_1(S(U_1,V_1))=T(U_1,V_1),
\end{equation}
where $\Psi_1$ acts coordinatewise.
The whole path is constant on $[1/8,1/2]$, so its completed right-continuous natural filtration satisfies
\begin{equation}\label{eq:example-filtration-pause}
 \cF_t^X=\cF_{1/4}^X\qquad(1/4\leq t\leq3/8).
\end{equation}
Moreover, $U_1,V_1$ are independent standard Gaussians, and $V_1$ is independent of $\cF_{1/4}^X$.
The terminal coordinate projection $x\mapsto x_1^Z\in\R^q$ is $1$-Lipschitz.
On the disjoint events $E_{k,j}:=\{\theta_{k,j}(U_1)=1\}\in\cF_{1/4}^X$, the Gaussian density on $[0,1]$ and \eqref{eq:counterexample-continuous-laws} give
\begin{equation}\label{eq:example-conditional-laws}
 \P(E_{k,j})\geq c\frac{|I_k|}{J_k},\qquad
 \cW_1\bigl(\Lc(Z_1\mid\cF_{1/4}^X),\mu_{k,j}\bigr)
 \leq r_k/100.
\end{equation}

\textbf{Step 4: Upper and lower bounds.}
Corollary~\ref{cor:osgood-vector-flow} verifies Assumptions~\ref{ass:osgood} and~\ref{ass:flow-rate} with $\beta=2/m$.
In particular, \eqref{eq:example-comparison-flow} gives $\Phi(a)\asymp[\log(1/a)]^{-2/m}$.
Theorem~\ref{thm:quantization}, with $p=1$, yields
\begin{equation}\label{eq:example-upper-bound}
 e^{\mathrm{ad}}_{M,1}(\Can(\Lc(X)))
 \leq C(\log M)^{-\eta_1},\qquad
 \eta_1=\min\left\{\frac12,\frac1{m+1}\right\}.
\end{equation}
Set $M_k:=\lfloor J_k/2\rfloor$.
By \eqref{eq:counterexample-law-separation}, \eqref{eq:example-filtration-pause}, and \eqref{eq:example-conditional-laws}, Corollary~\ref{cor:conditional-law-packing} applies with $t_0=1/4$, $\varepsilon_0=1/16$, radius $r_k/16$, and mass parameter $c|I_k|$.
For all sufficiently large $k$, it gives
\begin{equation}\label{eq:counterexample-lower-sequence}
 e^{\mathrm{ad}}_{M_k,1}(\Can(\Lc(X)))
 \geq c\,r_k2^{-k}.
\end{equation}

The choices in Steps~1--2 give $r_k2^{-k}=A2^{-k(q/m+2)}$ and $\log M_k\asymp2^{qk}$.
For $M_k\leq M<M_{k+1}$, monotonicity and \eqref{eq:counterexample-lower-sequence} give
\begin{equation}\label{eq:example-full-lower}
 e^{\mathrm{ad}}_{M,1}(\Can(\Lc(X)))
 \geq c_\beta2^{-(k+1)(q/m+2)}
 \geq c_\beta'(\log M)^{-1/m-2/q}.
\end{equation}
Since $q\geq m^2$ and $m=2/\beta$, we have $1/m+2/q\leq\beta/2+\beta^2/2$.
Combining \eqref{eq:example-upper-bound} and \eqref{eq:example-full-lower} proves \eqref{eq:counterexample-rate} in dimension $d=q+2=\lceil4/\beta^2\rceil+2$.
\end{proof}

\subsection{Proof of Proposition~\ref{prop:slow-quantization}}
\label{sec:slow-quantization-proof}

\begin{proof}[Proof of Proposition~\ref{prop:slow-quantization}]
We construct a two-dimensional diffusion $X=(U,Z)$ whose conditional future laws become increasingly complex in cells tending to spatial infinity.
The first coordinate selects a cell and later supplies its sampling randomness, while the second records a smooth path.
An interval of constant filtration separates selection from sampling.
We follow the steps used in the proof of Proposition~\ref{prop:counterexample}, replacing the Lipschitz seeds by a sampling mechanism within each cell in Step~2.

\textbf{Step 1: Separated laws and smooth samplers.}
Let $(n_k)_{k\geq k_0}$ be an increasing sequence of positive integers tending to infinity, to be chosen in Step~5.
Fix $\psi\in C_c^\infty((0,1))$ with $0\leq\psi\leq1$ and $\max\psi=1$.
At level $k\geq k_0$, divide $[3/4,1]$ into $n_k$ equal intervals and place a bump of height $r_k>0$ in each selected interval.
The resulting binary paths are
\begin{equation}
 f_{k,\omega}(t):=
 r_k\sum_{i=1}^{n_k}\omega_i
 \psi\bigl(4n_k(t-3/4)-(i-1)\bigr),\qquad
 \omega\in\{0,1\}^{n_k}.
\end{equation}
These $2^{n_k}$ smooth paths vanish on $[0,3/4]$ and satisfy
\begin{equation}\label{eq:slow-path-bounds}
 0\leq f_{k,\omega}\leq r_k,\qquad
 \|f_{k,\omega}'\|_\infty\leq4n_kr_k\|\psi'\|_\infty,\qquad
 \|f_{k,\omega}-f_{k,\omega'}\|_\infty=r_k
 \quad(\omega\ne\omega').
\end{equation}
We choose $r_k:=n_k^{-1}$ so that the time derivatives are uniformly bounded by $4\|\psi'\|_\infty$.
Lemma~\ref{lem:half-subset-packing} provides a constant $c_*>0$ such that, for sufficiently large $k_0$, there are at least $e^{c_*2^{n_k}}$ separated uniform laws on half-subsets of these paths.
We pick
$
 J_k:=\lfloor e^{c_*2^{n_k}}\rfloor
$
of these laws and denote them by $\mu_{k,j}\in\cP_1(\cC_1)$.
They satisfy
\begin{equation}\label{eq:slow-law-separation}
 \cW_1(\mu_{k,j},\mu_{k,j'})\geq r_k/4\qquad(j\ne j').
\end{equation}

Fix a tolerance $\epsilon_*\in(0,1/32)$.
Enumerate the $2^{n_k-1}$ atoms of $\mu_{k,j}$ as $f_{k,j,1},\ldots,f_{k,j,2^{n_k-1}}$.
In each interval $((\ell-1)/2^{n_k-1},\ell/2^{n_k-1})$ choose a smooth compactly supported cutoff $0\leq\theta_{k,\ell}\leq1$ equal to one on its central subinterval of relative length $1-\epsilon_*$, and put
\begin{equation}
 H_{k,j}(s,t):=\sum_{\ell=1}^{2^{n_k-1}}\theta_{k,\ell}(s)f_{k,j,\ell}(t).
\end{equation}
The disjoint supports imply that $H_{k,j}$ is smooth, vanishes for $s$ near $0$ and $1$, and satisfies
\begin{equation}\label{eq:slow-sampler-bounds}
 0\leq H_{k,j}\leq r_k,\qquad
 |\partial_tH_{k,j}|\leq4\|\psi'\|_\infty,\qquad
 H_{k,j}(s,t)=0\quad(t\leq3/4).
\end{equation}
For $R_0$ uniform on $(0,1)$, coupling to the path assigned to the interval containing $R_0$ gives
\begin{equation}\label{eq:slow-sampler-error}
 \cW_1\bigl(\Lc(H_{k,j}(R_0,\cdot)),\mu_{k,j}\bigr)
 \leq\epsilon_*r_k,
\end{equation}
since the paths agree on the cutoffs' plateaus and differ by at most $r_k$ elsewhere.

\textbf{Step 2: Sampling within a cell.}
Write $\Phi_{\mathrm G}$ for the standard Gaussian distribution function and choose a fixed $A>0$ such that
\begin{equation}\label{eq:slow-gaussian-shift}
 \frac{\operatorname{arctanh}(1/2)}{A\sqrt{2\pi}}\leq\epsilon_*.
\end{equation}
Starting at $u_0=c+\delta\tanh v_0$ in a cell $(c-\delta,c+\delta)$, consider
\begin{equation}
 \mathrm dU_t=-A^2h_1(t)^2\delta y_t(1-y_t^2)\,\mathrm dt
 +Ah_1(t)\delta(1-y_t^2)\,\mathrm dB_t,
 \qquad y_t:=\frac{U_t-c}{\delta},
\end{equation}
where $h_1$ is smooth, supported in $(1/2,5/8)$, and has squared integral one.
It\^o's formula applied to $\tanh$ gives
\begin{equation}\label{eq:slow-cell-sampling}
 U_t=c+\delta\tanh\left(v_0+A\int_{1/2}^t h_1(s)\,\mathrm dB_s\right),
 \qquad 1/2\leq t\leq5/8.
\end{equation}
This stays strictly inside the cell.
In particular, at the end of this stage,
\begin{equation}
 U=c+\delta\tanh(v_0+AG),\quad G\sim\mathcal N(0,1),\quad
 R:=\Phi_{\mathrm G}\left(\frac{\operatorname{arctanh}((U-c)/\delta)}A\right)
 =\Phi_{\mathrm G}(G+v_0/A).
\end{equation}
If $|u_0-c|\leq\delta/2$, then $|v_0|\leq\operatorname{arctanh}(1/2)$, and
\begin{equation}
 \TV\bigl(\Lc(R),\operatorname{Unif}(0,1)\bigr)
 \leq\TV\bigl(\mathcal N(v_0/A,1),\mathcal N(0,1)\bigr)
 \leq\frac{|v_0|}{A\sqrt{2\pi}}\leq\epsilon_*.
\end{equation}
Here $\TV(\mu,\nu)=\sup_D|\mu(D)-\nu(D)|$.

\textbf{Step 3: The two-dimensional diffusion.}
Assign the $J_k$ laws to equally spaced centers $c_{k,j}:=k+(j-1/2)/J_k$ in $(k,k+1)$.
Around each center, use a cell of radius $\delta_k:=1/(5J_k)$ and write $y_{k,j}(u):=(u-c_{k,j})/\delta_k$ for its normalized coordinate.
Since $\delta_k<1/(4J_k)$, the enlarged cells $(c_{k,j}-2\delta_k,c_{k,j}+2\delta_k)$ are disjoint and contained in $(k,k+1)$.
These cells have no finite accumulation point as $k\to\infty$.
Choose $\eta\in C_c^\infty((-2,2))$ with $0\leq\eta\leq1$ and $\eta=1$ on $[-1,1]$, and define
\begin{equation}
 a(u):=\sum_{k,j}\delta_k(1-y_{k,j}(u)^2)\eta(y_{k,j}(u)),\qquad
 d(u):=-\sum_{k,j}\delta_k y_{k,j}(u)(1-y_{k,j}(u)^2)\eta(y_{k,j}(u)).
\end{equation}
These sums are locally finite with at most one nonzero summand at each point, so $a,d$ are smooth and bounded.
For $|y_{k,j}(u)|<1$ put
\begin{equation}
 F_{k,j}(u,t):=H_{k,j}\left(
 \Phi_{\mathrm G}\left(\frac{\operatorname{arctanh}(y_{k,j}(u))}A\right),t\right),
\end{equation}
and extend it by zero outside that cell.
The vanishing of $H_{k,j}$ near its parameter endpoints makes this extension smooth.
Thus $F:=\sum_{k,j}F_{k,j}$ is smooth, $F(u,t)=0$ for $t\leq3/4$, and $|\partial_tF|\leq4\|\psi'\|_\infty$.

Choose smooth $h_0$ supported in $(0,1/8)$ with squared integral one.
Using $h_1$ from the cell sampling step, define $X=(U,Z)$ from $(0,0)$ by
\begin{equation}\label{eq:slow-two-dimensional-sde}
 \mathrm dU_t=A^2h_1(t)^2d(U_t)\,\mathrm dt
 +[h_0(t)+Ah_1(t)a(U_t)]\,\mathrm dB_t,\qquad
 \mathrm dZ_t=\partial_tF(U_t,t)\,\mathrm dt.
\end{equation}
The coefficients are smooth and bounded.
Local Lipschitz continuity and boundedness give a nonexplosive pathwise unique strong solution from every deterministic starting pair, verifying Assumption~\ref{ass:coefficients}.

\textbf{Step 4: The conditional laws and the lower bound.}
We have $U_{1/4}\sim\mathcal N(0,1)$.
The whole process is constant on $[1/8,1/2]$, so its completed right-continuous natural filtration satisfies
\begin{equation}\label{eq:slow-information-pause}
 \cF_t^X=\cF_{1/4}^X\qquad(1/4\leq t\leq3/8).
\end{equation}
For the disjoint events $E_{k,j}:=\{|U_{1/4}-c_{k,j}|\leq\delta_k/2\}\in\cF_{1/4}^X$, the Gaussian density gives the common mass bound
\begin{equation}\label{eq:slow-cell-mass}
 \P(E_{k,j})\geq\frac{w_k}{J_k},\qquad
 w_k:=\frac{e^{-(k+1)^2/2}}{5\sqrt{2\pi}}.
\end{equation}
On $E_{k,j}$, the cell sampling formula applies with $v_0=\operatorname{arctanh}(y_{k,j}(U_{1/4}))$.
Pathwise uniqueness identifies it with the first coordinate of \eqref{eq:slow-two-dimensional-sde}, which is constant after $5/8$.
Since $F$ vanishes up to $3/4$, integration of the second equation gives
\begin{equation}\label{eq:slow-output-path}
 Z_t=H_{k,j}(\Phi_{\mathrm G}(G+v_0/A),t),\qquad
 G:=\int_{1/2}^{5/8}h_1(s)\,\mathrm dB_s.
\end{equation}
The variable $G$ is standard Gaussian and independent of $\cF_{1/4}^X$.
The sampling estimate and \eqref{eq:slow-sampler-error}, together with the diameter bound $\|H_{k,j}(s,\cdot)-H_{k,j}(s',\cdot)\|_\infty\leq r_k$, give
\begin{equation}\label{eq:slow-conditional-approximation}
 \cW_1\bigl(\Lc(Z\mid\cF_{1/4}^X),\mu_{k,j}\bigr)
 \leq2\epsilon_*r_k\leq r_k/16\qquad\text{on }E_{k,j}.
\end{equation}
The coordinate projection $x\mapsto x^Z$ from $\cC_2$ to $\cC_1$ is $1$-Lipschitz.
Together with \eqref{eq:slow-law-separation}, \eqref{eq:slow-information-pause}, and \eqref{eq:slow-cell-mass}, this verifies all hypotheses of Corollary~\ref{cor:conditional-law-packing} with $t_0=1/4$, $\varepsilon_0=1/16$, packing radius $r_k/16$, and mass $w_k$.
Set $M_k:=\lfloor J_k/2\rfloor$.
Since $r_kw_k/32<\varepsilon_0$, Corollary~\ref{cor:conditional-law-packing} gives
\begin{equation}\label{eq:slow-sequence-lower}
 e^{\mathrm{ad}}_{M_k,1}(\Can(\Lc(X)))\geq\frac{r_kw_k}{32}.
\end{equation}

\textbf{Step 5: Choice of parameters.}
It remains to choose $n_k$.
By \eqref{eq:slow-cell-mass} and \eqref{eq:slow-sequence-lower}, we have
\begin{equation}
 e^{\mathrm{ad}}_{M_k,1}(\Can(\Lc(X)))
 \geq c\,n_k^{-1}e^{-(k+1)^2/2},\qquad
 \log M_k\asymp2^{n_k}.
\end{equation}
To make the growth of $\log M_k$ dominate the Gaussian level mass for every $\beta>0$, choose $n_k:=k^3$.
For $M_k\leq M<M_{k+1}$, monotonicity gives, for every $\beta>0$,
\begin{equation}
 (\log M)^\beta e^{\mathrm{ad}}_{M,1}(\Can(\Lc(X)))
 \geq c_\beta2^{\beta k^3}(k+1)^{-3}e^{-(k+2)^2/2}
 \longrightarrow\infty.
\end{equation}
This proves \eqref{eq:slow-quantization-rate}.
\end{proof}

\end{document}